\documentclass[12pt]{article}

\usepackage{a4wide,epsfig}
\usepackage{subfig}

\usepackage{amsfonts,amsmath}
\usepackage[utf8]{inputenc}   
\usepackage{float}
\usepackage{booktabs}
\usepackage[subpreambles=true]{standalone}
\newtheorem{theorem}{Theorem}[section]

\newtheorem{lemma}[theorem]{Lemma}

\newtheorem{remark}{Remark}[section]
\newcommand{\proof} [1]
   { \noindent {\bf Proof.} #1 \hfill\rule{0.5em}{1.2ex} \par\medskip}

\newcommand{\norm}[1]{{\left\lVert{#1}\right\rVert}}

\newcommand{\spf}[2]{{\left\langle{#1},{#2}\right\rangle}}

\newcommand{\pwc}{\mathcal{S}}
\newcommand{\OV}{\operatorname{\mathsf{V}}}
\newcommand{\OK}{\operatorname{\mathsf{K}}}

\numberwithin{equation}{section} 

\begin{document}

\setcounter{page}{1}

\title{Stable least-squares space-time boundary element methods for the wave
equation}
\author{Daniel~Hoonhout$^1$, Richard~L\"oscher$^2$, \\[1mm]
  Olaf~Steinbach$^2$, Carolina~Urz\'ua--Torres$^1$}
\date{$^1$Delft Institute of Applied Mathematics, TU Delft, \\[1mm]
  Mekelweg 4, 2628CD Delft, The Netherlands \\[2mm] 
  $^2$Institut f\"ur Angewandte Mathematik, TU Graz, \\[1mm]
  Steyrergasse 30, 8010 Graz, Austria}

\maketitle

\begin{abstract}
In this paper, we recast the variational formulation corresponding to the 
single layer boundary integral operator $\OV$ for the wave equation as 
a minimization problem in $L^2(\Sigma)$, where
$\Sigma := \partial \Omega \times (0,T)$ is the lateral boundary of the
space-time domain $Q := \Omega \times (0,T)$. For discretization, the
minimization problem is restated as a mixed saddle point formulation.
Unique solvability is established by combining conforming nested boundary
element spaces for the mixed formulation such that the related
bilinear form is discrete inf-sup stable. We 
analyze under which conditions the discrete inf-sup stability is satisfied,
and, moreover, we show that the mixed formulation provides a simple error
indicator, which can be used for adaptivity. We present several numerical
experiments showing the applicability of the method to different time-domain
boundary integral formulations used in the literature.
\end{abstract}

\section{Introduction}
Time-domain boundary integral equations and boundary element methods (BEM) for 
evolution problems are well established in the literature, see, for example 
\cite{Costabel:2017} for an overview. The common procedure is to 
either first discretize the temporal part using convolution quadrature, and 
then applying a boundary element method for the spatial variables, see,
e.g., \cite{AntesRubergSchanz:2004, BanjaiSauter:2008,
  Hassel:2017, Sayas:2016}; or to use BEM
with spatial basis functions and temporal basis functions chosen separately,
and then considered together as tensor product. We refer to
\cite{Aimi:2011, Aimi:2009, Aimi:2008, Gimperlein:2017, Gimperlein:2019,
  Gimperlein:2020, Gimperlein:2018, Ha-Duong:2003}, to name a few.
Typically, the choice of temporal basis function is done in order to obtain
a marching-on-in-time algorithm, which is an explicit time 
stepping scheme \cite{Terrasse:1993}.

Lately, the interest of discretizing space and time simultaneously has been 
increasing, resulting in so called space-time discretization methods. 
Admittedly, space-time discretizations lead to larger systems of algebraic 
equations to be solved. Nevertheless, these methods offer the advantage of 
having full control of the discretization in space and time at once, 
allowing for space-time adaptivity. In order to see this, it is worth 
noting that although space-time discretizations may lead to tensor-product 
basis functions on structured space-time meshes, they treat time as if it were 
another spatial variable, and thus, also permit unstructured meshes.
Moreover, space-time methods allow for 
preconditioning and parallelization in the space-time domain, which gives 
more flexibility in the construction of efficient solvers than time 
stepping methods, see, e.g., \cite{Gander:2015,GanderNeumueller:2016}.

The success of space-time BEM for parabolic problems
\cite{Dohr:2019,DohrNiinoSteinbach:2019,MertaOfWatschingerZapletal:2022} and
the promising developments for the wave equation when using space-time finite
element methods (FEM)
\cite{DoerflerFindeisenWieners:2016,LoescherSteinbachZank:2022,Perugia2018,
  SteinbachZank:2018,SteinbachZank:2021a}, encourage us to also study
space-time BEM for hyperbolic problems. For this, the first step is to
consider the literature on time-domain boundary integral equations for
the wave equation.

The standard approaches for BEM for the wave equation
started with the groundbreaking works of Bamberger and Ha-Duong
\cite{Bamberger:1986}, and Aimi et al.~\cite{Aimi:2009}. 
The main difficulty in the numerical analysis of these formulations is in 
the so-called norm gap, coming from continuity and coercivity estimates in 
different space-time Sobolev norms. When using the energetic BEM from 
\cite{Aimi:2009}, a complete stability and error analysis can be done in 
$L^2(\Sigma)$, see \cite{Joly:2017}, where $\Sigma := \partial \Omega
\times (0,T)$ is the lateral boundary of the space-time domain
$Q := \Omega \times (0,T)$. Hence, the energetic BEM is amenable 
to space-time discretizations. Its disadvantage, though, is that it requires 
the Dirichlet data to be sufficiently regular, i.e., in $H^1(\Sigma)$.

Using a generalized inf-sup stable variational formulation 
\cite{SteinbachZank:2021a} for the wave equation, in \cite{SteinbachUrzua:2021} 
we derived inf-sup stability conditions for all boundary integral operators
in related trace spaces, overcoming norm-gaps and also the need for extra 
regularity of the Dirichlet data. However, the standard discretization of the 
single layer boundary integral operator $\OV$ by means of space-time piecewise 
constant basis functions does not provide an inf-sup stable pair
\cite{Guardasoni:2010} in one spatial dimension,
which we believe will also be the case for $d=2,3$.

As an alternative, we proposed a regularisation via a modified Hilbert
transform \cite{SteinbachZank:2018}, the resulting composition with $\OV$
becomes elliptic in the natural energy space $[H^{1/2}_{,0}(\Sigma)]'$,
similarly to what is known for boundary integral operators for second-order
elliptic partial differential equations \cite{SteinbachUrzuaZank:2022}.
At the time of writing this article, this strategy had two drawbacks: 
it introduces further computational costs, and so far it is only applicable 
to space-time meshes that admit a tensor product structure. However, these
obstacles could be circumvented applying techniques used in \cite{SteinbachZank:2021}.

Another approach is to replace the straightforward variational formulation by
a least-squares/minimal residual equation. For FEM this
has been extensively studied for time independent problems, see Bochev and
Gunzberger \cite{BochevGunzburger:2006,BochevGunzburger:2016,
  BochevGunzburger:2009}. Time dependent parabolic problems have been
investigated in the context of first order least squares systems (FOSLS) in
\cite{FuehrerKarkulik:2021} and in the context of minimal residual
Petrov--Galerkin methods in \cite{Andreev:2013,StevensonWesterdiep:2021}.
In the context of BEM for elliptic partial differential equations this has
been studied in \cite{Steinbach:2023} and, recently, for FEM for the wave
equation in \cite{FuehrerGonzalesKarkulik:2023,
  KoetheLoescherSteinbach:2023}. In this paper we combine these ideas to also
have a least-squares boundary integral formulation that works for
hyperbolic problems. In addition to a stable method, we get an error
indicator that can be used for space-time adaptivity. It is worth pointing
out that, although we present the theory for the wave equation in one spatial
dimension, the underlying abstract framework is dimension-independent and, 
consequently, we expect the theory to carry over to higher dimensions.

This paper is organized as follows. Section~\ref{sec:LSVF} introduces our
notation and model problem. In particular, we remind the reader that in one
spatial dimension the single layer boundary integral operator $\OV$ is an
isomorphism from $[H^1_{,0}(\Sigma)]^{\prime}$ to $L^2(\Sigma)$, and we
derive the associated inf-sup constant. Then, we introduce a least-squares
variational formulation for the related boundary integral equation.
In Section~\ref{sec:MixedBEM}, we present the stable discretization of our
least squares formulation. For this, we propose a mixed BEM, and show its
unique solvability in Theorems~\ref{Thm Sh} and 
\ref{thm:discrete-stability-n-time-slices}. Moreover, Lemma~\ref{lem:aPriori} 
establishes the convergence of the method, and
Lemma~\ref{Lem ph error indicator} provides the conditions under which we
obtain a reliable error indicator. In Section~\ref{sec:NumExp}, we provide
numerical experiments to verify our theory. There, we investigate the
performance of the proposed least-squares formulation for three different
first-kind boundary integral equations related to our Dirichlet problem.
With this, we compare how different mapping and stability properties 
affect the numerical behaviour of the proposed method. We pay special attention 
to the requirements to have a reliable error indicator and the numerical study 
of the resulting adaptive scheme. Finally, we give some conclusions and
comment on ongoing work.

\section{Least-Squares Variational Formulation}
\label{sec:LSVF}

As in \cite{Aimi:2009,SteinbachUrzuaZank:2022}, we consider the Dirichlet 
boundary value problem for the homogeneous wave equation in the one-dimensional 
spatial domain $\Omega = (0,L)$, with zero initial conditions, and for a given
time horizon $T>0$,
\begin{equation}\label{DBVP}
  \left. \begin{array}{rcll}
           \partial_{tt} u(x,t) - \partial_{xx} u(x,t)
           & = & 0 & \mbox{for} \; (x,t) \in
                     Q := (0,L) \times (0,T), \\[1mm]
           u(x,0) = \partial_t u(x,t)_{|t=0}
           & = & 0 & \mbox{for} \; x \in (0,L),\\[1mm]
           u(0,t)
           & = & g_0(t) & \mbox{for} \; t \in (0,T), \\[1mm]
           u(L,t)
           & = & g_L(t) & \mbox{for} \; t \in (0,T).
         \end{array} \right \}
\end{equation}
In the one-dimensional case, the fundamental solution of the wave
equation is the Heaviside function
\[
U^*(x,t) = \frac{1}{2} \, {\mathsf{H}}(t-|x|),
\]
and we can represent the solution $u$ of (\ref{DBVP})
by using the single layer potential
\[
  u(x,t) := (\widetilde{\OV}w)(x,t) = \frac{1}{2} \int_0^{t-|x|} w_0(s) \, ds +
  \frac{1}{2} \int_0^{t-|x-L|} w_L(s) \, ds \quad \mbox{for} \;
  (x,t) \in Q .
\]
To determine the yet unknown density functions $(w_0,w_L)$, we consider
the boundary integral equations for $x \to 0$,
\begin{equation}\label{BIE 0}
(\OV_0w)(t) :=
  \frac{1}{2} \int_0^t w_0(s) \, ds +
\frac{1}{2} \int_0^{t-L} w_L(s) \, ds = g_0(t) \quad \mbox{for} \;
t \in (0,T),
\end{equation}
and for $ x \to L$,
\begin{equation}\label{BIE L}
  (\OV_Lw)(t) :=
\frac{1}{2} \int_0^{t-L} w_0(s) \, ds +
  \frac{1}{2} \int_0^t w_L(s) \, ds = g_L(t) \quad \mbox{for} \; t \in (0,T).
\end{equation}
We write the boundary integral equations (\ref{BIE 0}) and (\ref{BIE L}) in
compact form, for $w=(w_0,w_L)$, as
\begin{equation}\label{BIE}
  (\OV w)(t) =
  \begin{pmatrix}
           (\OV_0w)(t) \\
           (\OV_Lw)(t)
  \end{pmatrix}
  = \begin{pmatrix}
        \OV_{00} & \OV_{0L} \\ 
        \OV_{L0} & \OV_{LL}
    \end{pmatrix}
    \begin{pmatrix}
        w_0 \\ 
        w_L
    \end{pmatrix}
    (t)=
    \begin{pmatrix}
        g_0(t) \\ 
        g_L(t)
    \end{pmatrix}
    = g(t), \quad t \in (0,T) .
\end{equation}                         
In energetic BEM \cite{Aimi:2009}, instead
of (\ref{BIE}), the time derivative of (\ref{BIE}) is considered,
\begin{equation}\label{BIE energetic}
\partial_t (\OV w)(t) = \partial_t g(t) \quad \mbox{for} \; t \in (0,T).
\end{equation}
Recall the ellipticity estimate
\cite[Theorem 2.1]{Aimi:2009}, see also
\cite[Theorem 2.1]{SteinbachUrzuaZank:2022},
\[
  \langle w , \partial_t \OV w \rangle_{L^2(\Sigma)} \geq
  c_S(n) \, \| w \|^2_{L^2(\Sigma)} \quad \mbox{for all} \;
  w = (w_0,w_L) \in L^2(\Sigma) := L^2(0,T) \times L^2(0,T),
\]
where
\begin{equation*}\label{ellipticity constant}
  c_S(n) := \sin^2 \frac{\pi}{2(n+1)},
\end{equation*}
and 
\[n := \min \Big \{ m \in {\mathbb{N}} : T \leq m L \Big \},\]
is the number of time slices $T_j := ((j-1)L,j L)$ for $j =1,\ldots,n$
when $ T=nL$. In the case $ T < nL$, we define the last time slice
as $T_n:=((n-1)L,T)$, while all the others remain unchanged.

Since $\partial_t \OV : L^2(\Sigma) \to L^2(\Sigma)$ is also bounded,
unique solvability of the boundary integral equation
\eqref{BIE energetic} follows. 
Let us introduce $H^1_{0,}(\Sigma) := H^1_{0,}(0,T) \times H^1_{0,}(0,T)$, 
and note that 
$z \in H^1_{0,}(0,T)$ covers the zero initial condition $z(0)=0$.
Moreover, for $u = (u_0,u_L) \in H^1_{0,}(\Sigma)$ we have the norm definition
\[
  \| u \|^2_{H^1_{0,}(\Sigma)} := \| \partial_t u_0 \|^2_{L^2(0,T)} +
  \| \partial_t u_L \|^2_{L^2(0,T)} .
\]
Given that the time derivative $\partial_t : H^1_{0,}(\Sigma) \to L^2(\Sigma)$
is an isomorphism, e.g., \cite[Sect.~2.1]{SteinbachZank:2018}, we also 
have that $\OV : L^2(\Sigma) \to H^1_{0,}(\Sigma)$ is an isomorphism.

We define $H^1_{,0}(\Sigma) := H^1_{,0}(0,T) \times H^1_{,0}(0,T)$ in a
similar way, but with a zero terminal condition at $t=T$.
As in \cite[eqn.~(2.9)]{SteinbachUrzuaZank:2022} we also have the
ellipticity estimate
\[
  - \langle \overline{\partial}_t^{-1} \OV w,w \rangle_\Sigma \geq
  c_S(n) \, \| w \|^2_{[H^1_{,0}(\Sigma)]'} \quad \mbox{for all} \;
  w \in [H^1_{,0}(\Sigma)]',
\]
where
\[
  (\overline{\partial}_t^{-1}f)(t) = - \int_t^T f(s) \, ds, \quad
  t \in (0,T),
\]
is the inverse of $\partial_t : H^1_{,0}(\Sigma) \to L^2(\Sigma)$, and
$[H_{,0}^1(\Sigma)]'$ denotes the dual space of $H_{,0}^1(\Sigma)$ with
respect to $L^2(\Sigma)$, which is equipped with norm
\[
  \| w \|_{[H_{,0}^1(\Sigma)]'} = \sup_{0\neq v\in H_{,0}^1(\Sigma)}
  \frac{\langle w , v \rangle_\Sigma}{\| v \|_{H_{,0}^1(\Sigma)}},
\]
where $\langle \cdot , \cdot \rangle_\Sigma :
[H_{,0}^1(\Sigma)]' \times H_{,0}^1(\Sigma) \to {\mathbb{R}}$ denotes the
duality pairing as extension of the inner product
$\langle \cdot , \cdot \rangle_{L^2(\Sigma)}$ in $L^2(\Sigma)$.

Hence we conclude that $\OV : [H^1_{,0}(\Sigma)]' \to L^2(\Sigma)$ is an isomorphism and, 
in particular, bounded and satisfying the inf-sup stability condition. Our next aim is to find the 
inf-sup constant $\widetilde{c}_S(n)>0$ of
\begin{align*}
    \widetilde{c}_S(n)\norm{w}_{[H_{,0}^1(\Sigma)]'} = \sup_{0\neq q\in L^2(\Sigma)}
    \frac{\langle \OV w,q\rangle _{L^2(\Sigma)}}{\norm{q}_{L^2(\Sigma)}}
    =\norm{\OV w}_{L^2(\Sigma)}.
\end{align*} 
With this goal in mind, we first consider two auxiliary lemmas, where we
follow the ideas and notation of \cite{SteinbachZank:2018}.

\begin{lemma}\label{Lem representation w}
 Let $\{\mathcal{W}_k(t)\}_{k=0}^\infty:= \{ \cos(\alpha_k t)\}_{k=0}^\infty$ 
 with $\alpha_k:= \left(\frac{\pi}{2}+k\pi\right)\frac{1}{T}$, 
 and $w \in [H^1_{,0}(0,T)]'$. 
 Then there exists a unique $v\in L^2(0,T)$ such that $\partial_t v = w$ and
 \begin{equation*}\label{eq: representation w}
    \|w\|^2_{[H^1_{,0}(0,T)]'} = \frac{2}{T}\sum_{k=0}^\infty \alpha_k^{-2}
    \overline w_k^2 = \|v\|_{L^2(0,T)}^2,
 \end{equation*}
 where $\overline w_k := \langle w, \mathcal{W}_k \rangle_{(0,T)}$.
\end{lemma}
\proof{First, let us remind the reader that $\{\mathcal{W}_k(t)\}_{k=0}^\infty$ 
forms an orthogonal basis for $L^2(0,T)$ and $H^1_{,0}(0,T)$, while 
$\{\mathcal{V}_k(t)\}_{k=0}^\infty :=\{ \sin(\alpha_k t)\}_{k=0}^\infty$  
is an orthogonal basis for $L^2(0,T)$ and $H^1_{0,}(0,T)$.
Since $\partial_t: L^2(0,T) \to [H^1_{,0}(0,T)]'$ is an isomorphism, we can 
write $w = \partial_t v$ for a unique $v\in L^2(0,T)$. By representing $v$ 
by the orthogonal basis $\{\mathcal{V}_k\}_{k=0}^\infty$, i.e.,
$v = \sum_{k=0}^\infty v_k \mathcal{V}_k$, with
$v_k:= \frac{2}{T}\langle v, \mathcal{V}_k \rangle_{L^2(0,T)}$, we obtain
\begin{equation*}
  \norm{v}_{L^2(0,T)} = \left(\frac{T}{2}\sum_{k=0}^\infty |v_k|^2\right)^{1/2} ,
\end{equation*}
and further
\begin{equation}\label{w as dt_v}
  w = \partial_t v = \sum_{k=0}^\infty v_k \alpha_k \mathcal{W}_k.
\end{equation}
Using, that we can represent any function $q\in H_{,0}^1(0,T)$ as
$q(t) = \sum_{k=0}^\infty \overline{q}_k \mathcal{W}_k(t)$, with
$\overline{q}_k = \frac{2}{T}\spf{q}{\mathcal{W}_k}_{L^2(0,T)}$, we first
compute that 
\begin{align*}
  \norm{\partial_t q}_{L^2(0,T)} = \left( \frac{T}{2}
  \sum_{k=0}^\infty |\overline{q}_k|^2\alpha_k^2 \right)^{1/2}. 
\end{align*}
Now, by the definition of $\|\cdot \|_{[H^1_{,0}(0,T)]'}$, we have that 
\begin{eqnarray*}
  \|w\|_{[H^1_{,0}(0,T)]'}
  & = & \sup_{0\neq q \in H^1_{,0}(0,T)}
        \frac{\displaystyle \langle w, \sum_{k=0}^\infty \overline{q}_k
        \mathcal{W}_k \rangle_{(0,T)}}{\|\partial_t q \|_{L^2(0,T)}} 
        = \sup_{0 \neq q \in H^1_{,0}(0,T)}
        \frac{\displaystyle \sum_{k=0}^\infty \overline q_k \langle w,
        \mathcal{W}_k \rangle_{(0,T)}}
        {\displaystyle \left(\frac{T}{2} \sum_{k=0}^\infty
        \alpha_k^2 |\overline{q}_k|^2\right)^\frac{1}{2}} \\
  &=& \sup_{0 \neq q \in H^1_{,0}(0,T)}
      \frac{\displaystyle \sum_{k=0}^\infty \alpha_k\overline q_k
      \alpha_k^{-1}\overline w_k}
      {\displaystyle \left(\frac{T}{2} \sum_{k=0}^\infty \alpha_k^2
      |\overline{q}_k|^2\right)^\frac{1}{2}} 
      \leq \left(\frac{2}{T} \sum_{k=0}^\infty
      \overline w_k^2 \alpha_k^{-2}\right)^{1/2}.
\end{eqnarray*}
By picking
\[
  \hat q = \sum_{k=0}^\infty \alpha_k^{-2}\overline w_k \mathcal{W}_k
  \in H^1_{,0}(0,T),
\]
we can bound $\|w\|_{[H^1_{,0}(0,T)]'}$ from below by the same estimate, i.e., 
\begin{equation*}\label{w to v norm below}
  \|w\|_{[H^1_{,0}(0,T)]'} = \sup_{0\neq q \in H_{,0}^1(0,T)}
  \frac{\spf{w}{q}_{(0,T)}}{\norm{\partial_t q}_{L^2(0,T)}} \geq
  \frac{\spf{w}{\hat q}_{(0,T)}}{\|\partial_t \hat q\|_{L^2(0,T)}} =
  \left(\frac{2}{T}\sum_{k=0}^\infty \overline w_k^2 \alpha_k^{-2}\right)^{1/2}.
\end{equation*}
Thus, we have that
\begin{align*}
  \norm{w}_{[H_{,0}^1(0,T)]'} = \sqrt{\frac{2}{T}
  \left( \sum_{k=0}^\infty \overline w_k^2 \alpha_k^{-2}\right)}.
\end{align*}
By \eqref{w as dt_v}, we get $\overline{w}_k = \alpha_kv_k$ and compute,
\begin{equation*}\label{w to v final}
  \|w\|^2_{[H^1_{,0}(0,T)]'} =
  \frac{2}{T}\sum_{k=0}^\infty \alpha_k^{-2}\overline w_k^2 =
  \frac{T}{2} \sum_{k=0}^\infty |v_k|^2 = \|v\|^2_{L^2(0,T)}.
\end{equation*}
}

\begin{remark}
 The results of Lemma \ref{Lem representation w} also hold when considering the 
 lateral boundary $\Sigma$ as the domain, instead of $(0,T)$.
\end{remark}

\noindent
Using the compact form \eqref{BIE} for $w = (w_0,w_L)$, we define the operators 
$\OV_D $ and $\OV_{OD} $ as
\begin{equation}\label{eq:VD}
  \OV_D (t) := \begin{pmatrix} \OV_{00} & 0 \\ 0 & \OV_{LL} \end{pmatrix}
  \qquad \text{and } \qquad
  \OV_{OD} (t) := \begin{pmatrix} 0& \OV_{0L} \\ \OV_{L0} &0 \end{pmatrix}.
\end{equation}
Given Lemma \ref{Lem representation w} and definition \eqref{eq:VD}, we
can relate $\|\OV_D  w\|_{L^2(\Sigma)}$ to $\|w\|_{[H^1_{,0}(\Sigma)]'}$ as
summarised in the following lemma:
\begin{lemma}\label{Lem VD and w}
  Let $w \in [H^1_{,0}(\Sigma)]'$ and $\OV_D $ be defined as in \eqref{eq:VD}.
  Then
  \begin{equation*}\label{eq VD to w}
    \|\OV_D  w\|^2_{L^2(\Sigma)} = \frac{1}{4} \, \|w\|^2_{[H^1_{,0}(\Sigma)]'}.
  \end{equation*}
\end{lemma}
\proof{Let $w= \partial_t v$ for some $v \in L^2(\Sigma)$. By definition
  of $\OV_D $ we have
  \begin{equation*}\label{expression VD}
    \OV_D  w = \OV_D  \partial_t v = \frac{1}{2} v(t).
  \end{equation*}
  Hence, by Lemma \ref{Lem representation w}
  \begin{equation*}\label{Vd to w final}
    \|\OV_D w\|^2_{L^2(\Sigma)} =
    \frac{1}{4} \|v\|^2_{L^2(\Sigma)} = \frac{1}{4}\|w\|_{[H^1_{,0}(\Sigma)]'}^2.
  \end{equation*}
}

\noindent
Now we have all the tools to return to our study of the inf-sup 
constant $\tilde{c}_S(n)$ and provide the main result of this section.
\begin{theorem}\label{Thm ellipticity V*V}
  Let $w\in [H^1_{,0}(\Sigma)]'$ and let $n \in \mathbb{N}$ denote the number
  of time-slices. Then the operator $\OV:[H^1_{,0}(\Sigma)]' \to L^2(\Sigma)$
  is continuously bounded from below by constant $\widetilde{c}_S(n)$, i.e.,
  \begin{equation*}
    \|\OV w\|^2_{L^2(\Sigma)} \geq \widetilde{c}_S(n)^2
    \|w\|^2_{[H^1_{,0}(\Sigma)]'},
    \quad
    \mbox{where} \quad \widetilde{c}_S(n) :=
    \sin\left(\dfrac{\pi}{2(2n+1)}\right).
  \end{equation*}
\end{theorem}
\proof{We have
  \begin{eqnarray}\label{norm Vw decoupled}
    \| \OV w\|_{L^2(\Sigma)}
    &=& \langle (\OV_D  + \OV_{OD})w, (\OV_D  + \OV_{OD})w
        \rangle_{L^2(\Sigma)}\nonumber\\
    &=& \|\OV_D  w\|^2_{L^2(\Sigma)} + 2 \langle \OV_D  w, \OV_{OD}w
        \rangle_{L^2(\Sigma)} +  \|\OV_{OD}  w\|^2_{L^2(\Sigma)}.
  \end{eqnarray}
  Let $\Sigma_j$, $j=1,\ldots,n$, denote the lateral trace, restricted to
  the $j^{th}$ time-slice in time. By using the definitions \eqref{eq:VD},
  one can verify the following relation for all
  $w \in [H^1_{,0}(\Sigma)]'$:
  \begin{equation}\label{VD and VOD}
    \| \OV_D  w \|_{L^2(\Sigma_{j-1})} = \|\OV_{OD}  w\|_{L^2(\Sigma_{j})},
    \quad j=2,\ldots n.
  \end{equation}
  Using \eqref{norm Vw decoupled} and \eqref{VD and VOD}, we get
  \begin{equation*}
    \|\OV w\|^2_{L^2(\Sigma)}\geq \sum_{i=1}^n\|\OV_D  w\|^2_{L^2(\Sigma_i)} -
    2 \sum_{i=2}^{n} \|\OV_D  w\|_{L^2(\Sigma_i)}\|\OV_{D}w\|_{L^2(\Sigma_{i-1})}
    +\sum_{i=2}^{n}\|\OV_{D} w\|^2_{L^2(\Sigma_{i-1})},
  \end{equation*}
  which can be represented in matrix form as
  \begin{equation}\label{Vh norm matrix form}
    \|\OV w\|^2_{L^2(\Sigma)} \geq
    \left \langle
      \begin{pmatrix}
        2 & -1& & & \\
        -1 & 2 & -1 & &\\
        &\ddots &\ddots &\ddots & \\
        & & -1 & 2 & -1 \\
        & & & -1 & 1
      \end{pmatrix}
      \begin{pmatrix}
        \|\OV_D  w\|_{L^2(\Sigma_1)}\\
        \|\OV_D  w\|_{L^2(\Sigma_2)}\\
        \vdots \\
        \|\OV_D  w\|_{L^2(\Sigma_{n-1})}\\
        \|\OV_D  w\|_{L^2(\Sigma_{n})}
      \end{pmatrix},
      \begin{pmatrix}
        \|\OV_D  w\|_{L^2(\Sigma_1)}\\
        \|\OV_D  w\|_{L^2(\Sigma_2)}\\
        \vdots \\
        \|\OV_D  w\|_{L^2(\Sigma_{n-1})}\\
        \|\OV_D  w\|_{L^2(\Sigma_{n})}
      \end{pmatrix} \right\rangle.  
  \end{equation}
  The matrix in \eqref{Vh norm matrix form} corresponds to the one dimensional 
  finite difference matrix with (zero) initial Dirichlet condition and terminal 
  Neumann condition. The spectral properties of this matrix are henceforth well 
  understood, and its smallest eigenvalue is given by
  $\lambda_{min}= 4\sin^2 \left(\frac{\pi}{2(2n+1)}\right)$.
  Consequently, we obtain the bound
  \begin{eqnarray*}\label{Vw to VD eigenvalue}
    \|\OV w \|^2_{L^2(\Sigma)}
    &\geq& \lambda_{min} \sum_{i=1}^n\|\OV_D  w\|^2_{L^2(\Sigma_i)}\\
    &= &4\sin^2 \left(\frac{\pi}{2(2n+1)}\right) \|\OV_D w\|^2_{L^2(\Sigma)}
         \, \geq \,
         \sin^2 \left(\frac{\pi}{2(2n+1)}\right)\norm{w}_{[H_{,0}^1(\Sigma)]'}^2,
  \end{eqnarray*}
  where we applied Lemma \ref{Lem VD and w} in the last step.}

\noindent
As a direct consequence of Theorem \ref{Thm ellipticity V*V}, the inf-sup
stability condition is given by
\begin{equation}\label{inf sup V}
  \widetilde{c}_S(n) \, \| w \|_{[H^1_{,0}(\Sigma)]'} \leq
  \sup\limits_{0 \neq q \in L^2(\Sigma)}
  \frac{\langle \OV w , q \rangle_{L^2(\Sigma)}}{\| q \|_{L^2(\Sigma)}}
  \quad \mbox{for all} \; w \in H^1_{,0}(\Sigma),
\end{equation}
with the constant $\widetilde{c}_S(n)$ as given in
Theorem \ref{Thm ellipticity V*V}. 
In order to find its solution $w \in [H^1_{,0}(\Sigma)]'$, we consider
the minimization of the functional
\[
{\mathcal{J}}(v) := \frac{1}{2} \, \| \OV v - g \|^2_{L^2(\Sigma)},
\]
over all $v\in[H_{,0}^1(\Sigma)]$,
whose minimizer $w \in [H^1_{,0}(\Sigma)]'$ is determined as unique
solution of the gradient equation
\begin{equation}\label{gradient equation}
\OV^* \OV w = \OV^* g,
\end{equation}
where $\OV^* : L^2(\Sigma) \to H^1_{,0}(\Sigma)$ is the adjoint of
$\OV : [H^1_{,0}(\Sigma)]' \to L^2(\Sigma)$. When introducing the
adjoint $p := g - \OV w$, we end up with a mixed
variational formulation to find $p \in L^2(\Sigma)$ and
$w \in [H^1_{,0}(\Sigma)]'$ such that
\begin{equation}\label{mixed VF}
  \langle p , q \rangle_{L^2(\Sigma)} + \langle \OV w , q \rangle_{L^2(\Sigma)} =
  \langle g , q \rangle_{L^2(\Sigma)}, \quad
  \langle p , \OV v \rangle_{L^2(\Sigma)} = 0
\end{equation}
is satisfied for all $q \in L^2(\Sigma)$, and for all
$v \in [H^1_{,0}(\Sigma)]'$. In fact, the gradient equation
\eqref{gradient equation} is the Schur complement system of the
mixed formulation \eqref{mixed VF}. To establish unique solvability
of \eqref{gradient equation}, and therefore of \eqref{mixed VF},
we consider the Schur complement operator
$S := \OV^* \OV : [H^1_{,0}(\Sigma)]' \to H^1_{,0}(\Sigma)$.
For $w \in [H^1_{,0}(\Sigma)]'$ define $p_w = \OV w \in L^2(\Sigma)$ as unique
solution of the variational formulation
\begin{equation*}\label{Def pw}
  \langle p_w , q \rangle_{L^2(\Sigma)} =
  \langle \OV w , q \rangle_{L^2(\Sigma)} \quad \mbox{for all} \;
  q \in L^2(\Sigma) .
\end{equation*}
With this we obtain
\[
  \langle S w , w \rangle_\Sigma =
  \langle \OV^* \OV w , w \rangle_\Sigma =
  \langle p_w , \OV w \rangle_{L^2(\Sigma)} =
  \langle p_w, p_w \rangle_{L^2(\Sigma)} = \| p_w \|_{L^2(\Sigma)}^2,
\]
and the inf-sup stability condition \eqref{inf sup V} gives
\begin{equation}\label{bound w p}
  \widetilde{c}_S(n) \, \| w \|_{[H^1_{,0}(\Sigma)]'} \leq
  \sup\limits_{0 \neq q \in L^2(\Sigma)}
  \frac{\langle \OV w , q \rangle_{L^2(\Sigma)}}{\| q \|_{L^2(\Sigma)}} =
    \sup\limits_{0 \neq q \in L^2(\Sigma)}
  \frac{\langle p_w , q \rangle_{L^2(\Sigma)}}{\| q \|_{L^2(\Sigma)}}
  \leq \| p_w \|_{L^2(\Sigma)} .
\end{equation}
Hence, $S : [H^1_{,0}(\Sigma)]' \to H^1_{,0}(\Sigma)$ is elliptic
satisfying
\begin{equation*}\label{ellipticity S}
  \langle S w , w \rangle_\Sigma \geq [\widetilde{c}_S(n)]^2 \,
  \| w \|^2_{[H^1_{,0}(\Sigma)]'}
  \quad \mbox{for all} \; w \in [H^1_{,0}(\Sigma)]' .
\end{equation*}
From this, we conclude unique solvability of the gradient equation
\eqref{gradient equation} and of the mixed variational formulation
\eqref{mixed VF}.

\section{A Mixed Boundary Element Method}
\label{sec:MixedBEM}
Let
\[
  \pwc_H^0(\Sigma) := \pwc_{H,0}^0(0,T) \times \pwc_{H,L}^0(0,T) =
  \mbox{span} \{ \phi_\ell \}_{\ell=1}^{N_{H,0}} \times
  \mbox{span} \{ \phi_\ell \}_{\ell=N_{H,0}+1}^{N_H}
\]
and
\[
  \pwc_h^0(\Sigma) := \pwc_{h,0}^0(0,T) \times \pwc_{h,L}^0(0,T) =
  \mbox{span} \{ \varphi_k \}_{k=1}^{N_{h,0}} \times
  \mbox{span} \{ \varphi_k \}_{k=N_{h,0}+1}^{N_h}
\]
be two conforming nested boundary element spaces, i.e.,
\[
  \pwc_H^0(\Sigma) \subset \pwc_h^0(\Sigma) \subset
  L^2(\Sigma) \subset [H^1_{,0}(\Sigma)]',
\]
spanned by piecewise
constant basis functions $\phi_\ell$ and $\varphi_k$, which are defined
with respect to some nested decomposition of $\Sigma$ into boundary
elements $\tau_\ell^H$ and $\tau_k^h$ with local mesh sizes
$H_\ell$ and $h_k$, respectively. For $\tau_k^h \subset \tau_\ell^H$ we assume
$H_\ell = m h_k$ for some $m \in {\mathbb{N}}$. So we may define
a coarse grid mesh $\Sigma_H$ first, and any element
$\tau_\ell^H$ of $\Sigma_H$ is decomposed into $m$ equal sized elements
$\tau_k^h$ of the fine mesh $\Sigma_h$.

The Galerkin formulation of \eqref{mixed VF} is to find
$p_h \in \pwc_h^0(\Sigma)$ and $w_H \in \pwc_H^0(\Sigma)$ such that
\begin{equation}\label{mixed BEM}
  \langle p_h , q_h \rangle_{L^2(\Sigma)} +
  \langle \OV w_H , q_h \rangle_{L^2(\Sigma)} =
  \langle g , q_h \rangle_{L^2(\Sigma)}, \quad
  \langle p_h , \OV v_H \rangle_{L^2(\Sigma)} = 0
\end{equation}
is satisfied for all $q_h \in \pwc_h^0(\Sigma)$ and for all
$v_H \in \pwc_H^0(\Sigma)$. This is equivalent to a linear system of
algebraic equations,
\begin{equation}\label{LGS}
  \left( \begin{array}{cc}
           D_h & V_h \\
           V_h^\top &
         \end{array}
       \right)
       \left(
         \begin{array}{c}
           \underline{p} \\
           \underline{w}
         \end{array}
       \right)
       =
       \left(
         \begin{array}{c}
           \underline{g} \\
           \underline{0}
         \end{array}
         \right),
\end{equation}
where for $j,k=1,\ldots,N_h$ and
for $\ell=1,\ldots,N_H$ we have
\[
  D_h[j,k] = \int_\Sigma \varphi_k(x) \varphi_j(x) \, ds_x, \quad
  V_h[j,\ell] = \int_\Sigma (\OV\phi_\ell)(x) \varphi_j(x) \, ds_x, \quad
  g_j = \int_\Sigma g(x) \varphi_j(x) \, ds_x .
\]
Since the diagonal matrix $D_h$ is invertible, we can eliminate
$\underline{p}= D_h^{-1} [\underline{g}-V_h\underline{w}]$ to end up
with the Schur complement system
\begin{equation}\label{Schur complement system}
  S_h \underline{w} := V_h^\top D_h^{-1} V_h \underline{w} =
  V_h^\top D_h^{-1} \underline{g},
\end{equation}
which is nothing more than a Galerkin approximation of the gradient
equation \eqref{gradient equation}. By construction, the Schur
complement matrix $S_h = V_h^\top D_h^{-1} V_h$ is symmetric and at
least positive semi-definite. We will prove that the matrix $S_h$
is actually positive definite and hence that \eqref{Schur complement system}
and therefore \eqref{LGS} admits a unique solution.
\begin{theorem}\label{Thm Sh}
  Assume $T=L$, i.e., $n=1$. Let $\Sigma_H$ be a mesh of $\Sigma$, which 
  may be non-uniform and adaptive. Let $\Sigma_h$ be the
  fine mesh where each element $\tau_\ell^H$ of $\Sigma_H$ is decomposed
  into $m$ equal sized elements $\tau_k^h$. Then the Schur
  complement matrix $S_h^{L}$ is positive definite for all $m>2$, i.e.,
  \begin{equation*}\label{discrete ellipticity}
    (S_h^{L}\underline{w} , \underline{w}) \geq
    \Big( \frac{1}{2} - \frac{1}{m} \Big)^2 \,
    \| w_H \|_{[H^1_{,0}(\Sigma)]'}^2\quad \mbox{for all} \;
    \underline{w} \in {\mathbb{R}}^{N_H} \leftrightarrow
    w_H \in \pwc_H^0(\Sigma).
\end{equation*}
Here $S_h^{L}$ denotes the Schur complement matrix for a single time slice.
\end{theorem}
\proof{For $w_H \in \pwc_H^0(\Sigma) \subset [H^1_{,0}(\Sigma)]'$ the application
  of the Schur complement operator $S$ reads $Sw_H = \OV^*\OV w_H=\OV^* p_{w_H}$
  where $p_{w_H} \in L^2(\Sigma)$ is the unique solution of the variational
  formulation
  \[
    \langle p_{w_H} , q \rangle_{L^2(\Sigma)} =
    \langle \OV w_H , q \rangle_{L^2(\Sigma)} \quad \mbox{for all} \;
    q \in L^2(\Sigma).
  \]
  Now we consider the related Galerkin
  approximation $p_{w_H,h} \in \pwc_h^0(\Sigma)$
  as unique solution of the variational formulation
  \begin{equation}\label{Def ph}
    \langle p_{w_H,h} , q_h \rangle_{L^2(\Sigma)} =
    \langle p_{w_H} , q_h \rangle_{L^2(\Sigma)} =
    \langle \OV w_H , q_h \rangle_{L^2(\Sigma)} \quad \mbox{for all} \;
    q_h \in \pwc_h^0(\Sigma),
  \end{equation}
  i.e., we have to solve the linear system
  \[
    D_h \underline{p} = V_h \underline{w} \, .
  \]
  Instead of $Sw_H=\OV^*p_{w_H}$ we now define the approximation
  $\widetilde{S} w_H := \OV^* p_{w_H,h}$ for which we derive the
  matrix representation
  \[
S_h^L = V_h^\top D_h^{-1} V_h \, .
\]
Hence, we can write
\begin{eqnarray}\nonumber
  (S_h^L \underline{w} , \underline{w})
  & = & \langle \widetilde{S} w_H , w_H \rangle_{\Sigma} =
        \langle \OV^* p_{w_H,h} , w_H \rangle_{\Sigma} \\
  & = & \langle p_{w_H,h},\OV w_H \rangle_{L^2(\Sigma)} =
        \langle p_{w_H,h},p_{w_H,h} \rangle_{L^2(\Sigma)} =
        \| p_{w_H,h} \|^2_{L^2(\Sigma)} . \label{estimate Sh}
\end{eqnarray}
From the triangle inequality
\[
  \| p_{w_H} \|_{L^2(\Sigma)} =
  \| p_{w_H} - p_{w_H,h} + p_{w_H,h} \|_{L^2(\Sigma)} \leq
  \| p_{w_H} - p_{w_H,h} \|_{L^2(\Sigma)} + \| p_{w_H,h} \|_{L^2(\Sigma)}
\]
we get, by using \eqref{bound w p},
\begin{eqnarray}\nonumber
  \| p_{w_H,h} \|_{L^2(\Sigma)}
  & \geq & \| p_{w_H} \|_{L^2(\Sigma)} -
           \| p_{w_H} - p_{w_H,h} \|_{L^2(\Sigma)} \\
  & \geq & \widetilde{c}_S(n) \, \| w_H \|_{[H^1_{,0}(\Sigma)]'} -
           \| p_{w_H} - p_{w_H,h} \|_{L^2(\Sigma)}, \label{error p ph}
 \end{eqnarray}
 and it remains to estimate the approximation error of
 $ \| p_{w_H} - p_{w_H,h} \|_{L^2(\Sigma)}$.

 In the case $T=L$, the application of the wave single layer boundary
 integral operator $\OV$ is decoupled, i.e.,
 \[
   p_{w_H,0}(t) =  \frac{1}{2} \int_0^t w_{H,0}(s) \, ds, \quad
   p_{w_H,L}(t) =  \frac{1}{2} \int_0^t w_{H,L}(s) \, ds,
   \quad t \in (0,T).
 \]
 For the coefficients of the piecewise constant approximation at $x=0$,
 \[
  p_{w_H,0,h}(t) = \sum\limits_{k=1}^{N_{h,0}} p_{0,k} \varphi_k(t),
\]
we find from \eqref{Def ph} that
\[
  p_{0,k} = \frac{1}{h_k} \int_{t_{k-1}}^{t_k} p_{w_H,0}(s) \, ds \, 
  \quad \mbox{for} \; k=1,\ldots,N_{h,0}.
\]
By using standard arguments, see, e.g., \cite{Steinbach:2008}, and
$p_{w_H,0}'(t) = \frac{1}{2}w_{H,0}(t)$,
we obtain the error estimate
\[
  \int_{t_{k-1}}^{t_k} [p_{w_H,0}(t)-p_{w_H,0,h}(t)]^2 \, dt \leq
  \frac{1}{3} \, h_k^2
  \int_{t_{k-1}}^{t_k} [p_{w_H,0}'(t)]^2 dt =
  \frac{1}{12} \, h_k^2
  \int_{t_{k-1}}^{t_k} [w_{H,0}(t)]^2 dt \, ,
\]
and summing up this gives
\[
  \int_0^T [p_{w_H,0}(t)-p_{w_H,0,h}(t)]^2 \, dt \leq \frac{1}{12}
  \sum\limits_{k=1}^{N_{h,0}} h_k^2 \int_{t_{k-1}}^{t_k} [w_{H,0}(t)]^2 dt \, .
\]
When inserting
\[
w_{H,0}(t) = \sum\limits_{\ell=1}^{N_{H,0}} w_\ell \phi_\ell(t),
\]
assembling all fine grid contributions from the elements
$\tau_k^h \subset \tau_\ell^H$, and using $H_\ell = m h_k$, we further
conclude
\[
  \int_0^T [p_{w_H,0}(t)-p_{w_H,0,h}(t)]^2 \, dt \leq
  \frac{1}{12} \, \frac{1}{m^2}
  \sum\limits_{\ell=1}^{N_{H,0}} H_\ell^3 w_\ell^2 \, .
\]
By doing the same computations at $x=L$, and summing up both
contributions, this gives
\[
  \| p_{w_H} - p_{w_H,h} \|^2_{L^2(\Sigma)} \leq
  \frac{1}{12} \, \frac{1}{m^2} \, \sum\limits_{\ell=1}^{N_H}
  H_\ell^3 w_\ell^2 \, .
\]
We now consider a piecewise quadratic function
\[
v_H(t) = \sum\limits_{\ell=1}^{N_H} w_\ell \psi_\ell(t),
\]
where the bubble function $\psi_\ell$ in the boundary element
$\tau_\ell^H$ is defined by its form function
\[
\psi(s) = s (H-s) \quad \mbox{for} \; s \in (0,H).
\]
For this we compute
\[
\int_{\tau_\ell^H} \psi_\ell(t) \, dt =
  \frac{1}{6} \, H_\ell^3, \qquad
  \int_{\tau_\ell^H} [\psi_\ell'(t)]^2 \, dt
   = \frac{1}{3} \, H_\ell^3 .
\]
Thus, we have
\[
  \langle w_H , v_H \rangle_{L^2(\Sigma)} =
  \frac{1}{6} \, \sum\limits_{\ell=1}^{N_H} w_\ell^2 \, H_\ell^3,
\quad \text{ and } \quad
  \| v_H' \|^2_{L^2(\Sigma)} =
  \frac{1}{3} \, \sum\limits_{\ell=1}^{N_H} w_\ell^2 \, H_\ell^3 .
\]
With this, we finally obtain
\begin{eqnarray*}
  \| p_{w_H} - p_{w_H,h} \|_{L^2(\Sigma)}
  & \leq & \frac{1}{m} \sqrt{
           \frac{1}{12} \sum\limits_{\ell=1}^{N_H} H_\ell^3 w_\ell^2} =
           \frac{1}{m}
           \frac{\frac{1}{12} \sum\limits_{\ell=1}^{N_H} H_\ell^3 w_\ell^2}
           {\sqrt{\frac{1}{12} \sum\limits_{\ell=1}^{N_H} H_\ell^3 w_\ell^2}} 
           = \frac{1}{m} \, 
           \frac{\frac{1}{2} \, \langle w_H,v_H \rangle_{L^2(\Sigma)}}
           {\frac{1}{2} \, \| v_H' \|_{L^2(\Sigma)}} \\
  & \leq & \frac{1}{m} \,  
           \sup\limits_{0 \neq v \in H^1_{,0}(\Sigma)}
        \frac{\langle w_H,v \rangle_{L^2(\Sigma)}}{\| v'\|_{L^2(\Sigma)}} =
        \frac{1}{m} \, \| w_H \|_{[H^1_{,0}(\Sigma)]'} \, .
\end{eqnarray*}
Hence, we can write \eqref{error p ph} as
\begin{equation}\label{estimate pwh wh}
  \| p_{w_H,h} \|_{L^2(\Sigma)}
  \, \geq \, \Big( \widetilde{c}_S(n) - \frac{1}{m} \Big) \,
  \| w_H \|_{[H^1_{,0}(\Sigma)]'},
\end{equation}
which is positive for
\[
  \frac{1}{m} < \widetilde{c}_S(n) =
  \sin\left(\frac{\pi}{2(2n+1)}\right) \stackrel{n=1}{=}
  \sin \frac{\pi}{6} = \frac{1}{2}, \quad \mbox{i.e., for} \; m > 2.
\]
The assertion now follows from \eqref{estimate Sh}.}

\noindent 
Theorem \ref{Thm Sh} ensures unique solvability of the linear system
\eqref{Schur complement system}, and therefore of the mixed variational
formulation \eqref{mixed BEM} in the particular case $T=L$. But this
result can be generalized as follows.
\begin{theorem}\label{thm:discrete-stability-n-time-slices}
  Let $T=nL$ for some $n \in {\mathbb{N}}$ that induces time slices
  $((j-1)L,jL)$ for $j=1,\ldots,n$.  Let $\Sigma_H$ be a 
  uniform mesh of $\Sigma$. Let
  $\Sigma_h$ be the fine mesh where each element $\tau_\ell^H$ of
  $\Sigma_H$ is decomposed into $m$ equal sized elements $\tau_k^h$.
  We assume that $jL$, $j=0,\ldots,n$, are nodes of the 
  mesh $\Sigma_H$ at $x=0$ and at $x=L$, respectively. Then,
  \begin{equation*}\label{discrete ellipticity n}
    (S_h^{nL}\underline{w} , \underline{w}) \geq
    4\sin^2\left(\frac{\pi}{2(2n +1)}\right) \,
    \Big( \frac{1}{2} - \frac{1}{m} \Big)^2 \,
    \| w_H \|_{[H^1_{,0}(\Sigma)]'}^2\quad \mbox{for all} \;
    \underline{w} \in {\mathbb{R}}^{N_H} \leftrightarrow
    w_H \in \pwc_H^0(\Sigma).
  \end{equation*}
  Here, $S_h^{nL}$ denotes the Schur complement matrix for $n$ time slices.
\end{theorem}
\proof{Let $Q_h$ denote the $L^2$ projection with respect to the fine
    mesh $\Sigma_h$, defined as
    \begin{equation}\label{def:L2projection}
      \langle Q_h u , v_h \rangle_{L^2(\Sigma)} =
      \langle u , v_h \rangle_{L^2(\Sigma)}, \quad
      \text{for all } v_h \in \pwc^0_h(\Sigma),
    \end{equation}
    when $u \in L^2(\Sigma)$ is given.
    In case of a uniform refinement, we retain an equality analogous to
    \eqref{VD and VOD}, i.e.,
    \begin{equation}\label{QhVD to QhVOD}
      \|Q_h \OV_D  w\|_{L^2(\Sigma_{j-1})} =
      \|Q_h \OV_{OD}  w\|_{L^2(\Sigma_{j})}, \quad j=2,\dots,n.
    \end{equation}    
    Hence, we get
    \begin{align}\label{Sh to Qh V}
      (S_h^{nL}\underline{w} , \underline{w}) =
      (V_h^T D_h^{-1}V_h, \underline{w} , \underline{w})
      =
      \left(\sup_{0\neq q_h \in \pwc_h^0(\Sigma)}
      \frac{\langle \OV w_H, q_h \rangle_{L^2(\Sigma)}}{\|q_h\|_{L^2(\Sigma)}}
      \right)^2 = \|Q_h \OV  w_H\|_{L^2(\Sigma)}^2 .
    \end{align}
    Given \eqref{QhVD to QhVOD}, and following the lines of the proof of
    Theorem \ref{Thm ellipticity V*V}, we get
    \begin{equation}\label{QhV to QhVD}
      \|Q_h \OV  w_H\|_{L^2(\Sigma)}^2 \geq
      4 \, [\widetilde{c}_S(n)]^2 \, \|Q_h \OV_D  w_H\|^2_{L^2(\Sigma)}.
    \end{equation}
    Next, we make the observation that $ \|Q_h \OV_D  w_H\|^2_{L^2(\Sigma)}$
    is equivalent to considering the discretized operator $Q_h \OV$ on one
    time-slice, which suggests that we can apply Theorem \ref{Thm Sh}.
    The result now follows from this observation, \eqref{Sh to Qh V} and
    \eqref{QhV to QhVD}, 
    \begin{align}\label{eq:SnL-ellipticity}
      (S_h^{nL}\underline{w} , \underline{w})
      & = \|Q_h \OV  w_H\|_{L^2(\Sigma)}^2 
        \geq \, 4 \, [\widetilde{c}_S(n)]^2 \,
        \|Q_h \OV_D  w_H\|^2_{L^2(\Sigma)} \nonumber\\
      & \geq \, 4 \, [\widetilde{c}_S(n)]^2 \
        \Big( \frac{1}{2} - \frac{1}{m} \Big)^2 \,
        \| w_H \|_{[H^1_{,0}(\Sigma)]'}^2.
\end{align}
}

\begin{remark}\label{rem:meshconstraints}
  The relation \eqref{QhVD to QhVOD} does not need to hold on non-uniform 
  refinements. However, in absence of rounding errors, adaptively refined
  meshes should retain the property \eqref{QhVD to QhVOD} if one starts
  with a uniform initial mesh, and, as a consequence, the mixed boundary
  element method remains discrete inf-sup stable. We will see later, in
  the numerical experiments, how numerical errors do break this condition
  and how to remedy it.
\end{remark}

\begin{remark}
  Note, that in the limit case $h\to 0$, we have $m\to \infty$ and the
  bound of Theorem \ref{thm:discrete-stability-n-time-slices} becomes
  exactly the bound of Theorem \ref{Thm ellipticity V*V} in the
  continuous case.    
\end{remark}

\noindent
It remains to provide an a priori error estimate for the
unique solution of the mixed variational formulation \eqref{mixed BEM}.
Although this follows as in the elliptic case for the Laplace
equation \cite{Steinbach:2023}, here we present the main steps:
\begin{lemma}
\label{lem:aPriori}
  Let the assumption of Theorem \ref{thm:discrete-stability-n-time-slices} hold. Then, for the unique solutions $w\in [H_{,0}^1(\Sigma)]'$ of \eqref{BIE} and
  $w_H\in \pwc_H^0(\Sigma)$ of \eqref{mixed BEM}, there holds
\begin{align*}
	\norm{w-w_H}_{[H_{,0}^1(\Sigma)]'}\leq \left(1+\frac{2mc_2^{\OV}}{\sin\left(\frac{\pi}{2(2n+1)}\right)(m-2)}\right)\inf_{v_H\in \pwc_H^0(\Sigma)}\norm{w-v_H}_{[H_{,0}^1(\Sigma)]'}. 
\end{align*}  
\end{lemma}
\proof{When combining \eqref{Sh to Qh V} and \eqref{eq:SnL-ellipticity} we immediately obtain the discrete inf-sup stability condition
  \[
    2\widetilde{c}_S(n)\left( \frac{1}{2} - \frac{1}{m} \right) \,
    \| v_H \|_{[H^1_{,0}(\Sigma)]'} \leq
    \sup\limits_{0 \neq q_h \in \pwc_h^0(\Sigma)}
    \frac{\langle \OV v_H , q_h \rangle_{L^2(\Sigma)}}
    {\| q_h \|_{L^2(\Sigma)}} \quad \mbox{for all} \; v_H \in \pwc_H^0(\Sigma).
  \]
  Then, for the solution $w_H$ of \eqref{mixed BEM} and for
  arbitrary $v_H\in \pwc_H^0(\Sigma)$, we obtain, by using the triangle inequality,
  \[
    \| w-w_H \|_{[H^1_{,0}(\Sigma)]'} \leq
    \| w-v_H \|_{[H^1_{,0}(\Sigma)]'} +
    \| v_H-w_H \|_{[H^1_{,0}(\Sigma)]'}.
  \]
  Now, for the second term, we can use the 
  discrete inf-sup stability condition and \eqref{mixed BEM} for
  $g=\OV w$ to estimate
  \begin{eqnarray*}
    2\widetilde{c}_S(n)\left( \frac{1}{2} - \frac{1}{m} \right) \,
    \| v_H-w_H \|_{[H_{,0}^1(\Sigma)]'}
    & \leq & \sup_{0\neq q_h\in \pwc_h^0(\Sigma)}
             \frac{\langle \OV(v_H-w_H),q_h \rangle_{L^2(\Sigma)}}
             {\| q_h \|_{L^2(\Sigma)}} \\
    &&  \hspace*{-5cm} = \sup_{0\neq q_h\in \pwc_h^0(\Sigma)}
          \frac{\langle \OV v_H - (g-p_h),q_h \rangle_{L^2(\Sigma)}}
          {\| q_h \|_{L^2(\Sigma)}} \, = \,
       \sup_{0\neq q_h \in \pwc_h^0(\Sigma)}
          \frac{\langle \OV(v_H-w)+p_h , q_h \rangle_{L^2(\Sigma)}}
          {\| q_h \|_{L^2(\Sigma)}} \\[1mm]
    && \hspace*{-5cm} \leq \, c_2^{\OV} \, \| v_H-w \|_{[H_{,0}^1(\Sigma)]'} +
             \| p_h \| _{L^2(\Sigma)}.           
  \end{eqnarray*}
  Thus, it remains to estimate $\| p_h \|_{L^2(\Sigma)}$. Therefore, we
  consider \eqref{mixed BEM} with $q_h=p_h$, to get
  \begin{eqnarray*}
    \| p_h \|_{L^2(\Sigma)}^2
    & = & \langle p_h , p_h \rangle_{L^2(\Sigma)} \, = \,
          \langle g-Vw_H , p_h \rangle_{L^2(\Sigma)} \, = \,
          \langle \OV(w-w_H) , p_h \rangle_{L^2(\Sigma)} \\
    & = & \langle \OV(w-v_H),p_h \rangle_{L^2(\Sigma)} +
          \langle \OV(v_H-w_H),p_h \rangle_{L^2(\Sigma)} \\
    & = & \langle \OV(w-v_H),p_h \rangle_{L^2(\Sigma)} \, \leq \,
          c_2^{\OV} \, \| w-v_H \|_{[H_{,0}^1(\Sigma)]'}
          \| p_h \|_{L^2(\Sigma)}.
  \end{eqnarray*}
  Now, combining the estimates and taking the infimum over all
  $v_H\in \pwc_H^0(\Sigma)$ we obtain
  \begin{eqnarray*}
    \| w-w_H \|_{[H_{,0}^1(\Sigma)]'}
    & \leq & \left( 1 + \frac{c_2^{\OV}}{\widetilde{c}_S(n)\left(\frac{1}{2}-\frac{1}{m}\right)}
             \right) \inf_{v_H \in \pwc_H^0(\Sigma)}
             \| w-v_H \|_{[H^1_{,0}(\Sigma)]'}\\
    & = & \left( 1 + \frac{2mc_2^{\OV}}{\widetilde{c}_S(n)(m-2)} \right)
          \inf_{v_H\in \pwc_H^0(\Sigma)} \| w-v_H \|_{[H^1_{,0}(\Sigma)]'}. 
  \end{eqnarray*}}

\newpage
\noindent
Let $0<h<\underline {H}<H$ be given such that the inclusion
$\pwc_H^0(\Sigma) \subset \pwc_{\underline H}^0(\Sigma)\subset \pwc_h^0(\Sigma)$
holds, and assume that there exists $\widetilde{c}_S(n)>0$ such that
\begin{align}\label{eq:inf-sup-stability-Hu}
  \widetilde{c}_S(n) \, \| v_{\underline H} \|_{[H^1_{,0}(\Sigma)]'}
  \leq \sup_{0 \neq q_h\in \pwc_h^0(\Sigma)}
  \frac{\langle \OV v_{\underline H} , q_h \rangle_{L^2(\Sigma)}}
  {\| q_h \|_{L^2(\Sigma)}} \quad
  \text{for all }v_{\underline H}\in \pwc_{\underline H}^0(\Sigma)
\end{align}
is satisfied. Then \eqref{mixed BEM} admits a unique solution
$(\overline p_h,w_{\underline H})\in \pwc_h^0(\Sigma) \times
\pwc_{\underline H}^0(\Sigma)$. Note that, due to the inclusion
$\pwc_H^0(\Sigma)\subset \pwc_{\underline H}^0(\Sigma)$, we have that 
\begin{align*}
  \widetilde{c}_S(n) \, \| v_H \|_{[H^1_{,0}(\Sigma)]'}
  \leq \sup_{0\neq q_h\in \pwc_h^0(\Sigma)}
  \frac{\langle \OV  v_H,q_h \rangle_{L^2(\Sigma)}}{\| q_h \|_{L^2(\Sigma)}}
\end{align*} 
holds true for all $v_H\in \pwc_H^0(\Sigma)$, and thus \eqref{mixed BEM} also 
admits a unique solution $(p_h,w_{H})\in \pwc_h^0(\Sigma)\times
\pwc_{H}^0(\Sigma)$. Under a saturation assumption, we can now show that
$p_h\in \pwc_h^0(\Sigma)$ is an error estimator.

\begin{lemma}\label{Lem ph error indicator}
  Let $w\in [H_{,0}^1(\Sigma)]'$ be the unique solution of \eqref{BIE}.
  Further, let $(p_h,w_{\underline H})\in \pwc_h^0(\Sigma)\times
  \pwc_{\underline H}^0(\Sigma)$ and $(p_h,w_H)\in \pwc_h^0(\Sigma)\times
  \pwc_H^0(\Sigma)$ be the unique solution of 
  \eqref{mixed BEM}. If the saturation assumption 
  \begin{align}\label{eq:saturation-assumption}
    \| w-w_{\underline H} \|_{[H^1_{,0}(\Sigma)]'} \leq
    \delta \, \| w-w_H\|_{[H^1_{,0}(\Sigma)]},\quad \text{for }\delta \in (0,1)
  \end{align} 
  holds, then
  \begin{align*}
    \frac{1}{c_2^{\OV}} \, \| p_h \|_{L^2(\Sigma)} \leq
    \| w-w_H \|_{[H^1_{,0}(\Sigma)]}
    \leq
    \frac{c_2^{\OV}}{[\widetilde{c}_S(n)]^2} \, \frac{1}{1-\delta} \,
    \| p_h \|_{L^2(\Sigma)}.
  \end{align*} 
\end{lemma}

\proof{First, using \eqref{mixed BEM} and $\OV w=g$ and the boundedness
  of $\OV$, we compute
  \begin{eqnarray*}
    \| p_h \|^2_{L^2(\Sigma)}
    & = & \langle p_h , p_h \rangle_{L^2(\Sigma)} =
          \langle g-\OV w_H , p_h \rangle_{L^2(\Sigma)}\\
    & = & \langle \OV(w-w_H),p_h \rangle_{L^2(\Sigma)}
          \leq c_2^{\OV} \, \| w-w_H \|_{[H_{,0}^1(\Sigma)]} \,
          \| p_h \|_{L^2(\Sigma)},
  \end{eqnarray*}
  from which we conclude the first bound. To bound the error by $\norm{p_h}_{L^2(\Sigma)}$, let
  us first estimate
  \begin{align*}
    \| w-w_H \|_{[H^1_{,0}(\Sigma)]'} \leq
    \| w-w_{\underline H} \|_{[H^1_{,0}(\Sigma)]'} +
    \| w_{\underline H}-w_H \|_{[H^1_{,0}(\Sigma)]'}. 
  \end{align*}
  With the saturation \eqref{eq:saturation-assumption}, we conclude 
  \begin{align*}
    \| w-w_H \|_{[H^1_{,0}(\Sigma)]'} \leq \frac{1}{1-\delta} \,
    \| w_{\underline H}-w_H \|_{[H^1_{,0}(\Sigma)]'}.
  \end{align*}
  Thus, it is sufficient to bound the discrete error. We note that
  $w_{\underline H}-w_H\in \pwc_{\underline H}^0(\Sigma)$ and we can use the
  discrete inf-sup stability \eqref{eq:inf-sup-stability-Hu}, and
  together with \eqref{mixed BEM} we get
  \begin{eqnarray*}
    \widetilde{c}_S(n) \, \| w_{\underline H}-w_H \|_{[H^1_{,0}(\Sigma)]'}
    & \leq & \sup_{0\neq q_h\in \pwc_h^0(\Sigma)}
             \frac{\langle \OV(w_{\underline H}-w_H),q_h \rangle_{L^2(\Sigma)}}
             {\| q_h \|_{L^2(\Sigma)}}\\
    & = & \sup_{0\neq q_h\in \pwc_h^0(\Sigma)}
          \frac{\langle p_h-\overline p_h,q_h \rangle_{L^2(\Sigma)}}
          {\| q_h \|_{L^2(\Sigma)}} =
          \| p_h-\overline p_h \|_{L^2(\Sigma)}. 
  \end{eqnarray*}
  We can further bound this term by using again \eqref{mixed BEM} as follows
  \begin{eqnarray*}
    \| p_h-\overline p_h \|_{L^2(\Sigma)}^2
    & = & \langle p_h-\overline p_h,p_h-\overline p_h\rangle_{L^2(\Sigma)} =
          \langle \OV(w_{\underline H}-w_H),p_h-\overline p_h
          \rangle_{L^2(\Sigma)}\\[1mm]
    & = & \langle \OV(w_{\underline H}-w_H,p_h \rangle_{L^2(\Sigma)} \leq
          c_2^{\OV} \, \| w_{\underline H}-w_H \|_{[H^1_{,0}(\Sigma)]'}
          \| p_h \|_{L^2(\Sigma)} .
  \end{eqnarray*}
  Altogether, we now obtain that
  \begin{align*}
    \| w_{\underline H}-w_H \|_{[H^1_{,0}(\Sigma)]'} \leq
    \frac{c_2^{\OV}}{[\widetilde{c}_S(n)]^2} \, \| p_h \|_{L^2(\Sigma)},
  \end{align*}
  which concludes the proof. }

\begin{remark}\label{rem ph practice}
  The solution $(\overline p_h,w_{\underline H})\in \pwc_h^0(\Sigma) \times
  \pwc_{\underline H}^0(\Sigma)$ is only needed for the proof of the error
  estimator and does not need to be computed. In general, the idea is to
  have a stable method and then refine the mesh of the dual variable once
  more to get an error estimator. In particular, if the method is stable
  for the choices $\underline{H}=H/2$ then the choice
  $h=H/4$ gives an error estimator (or merely
  $\underline{H}=H$ is stable then $h=H/2$ gives
  an estimator). This is in some sense a generalization of the $h-h/2$ error
  estimator for elliptic equations. The behavior is also resembled by our
  numerical examples, as the choice $m=2$ gives a method that is stable
  in the primal variable, thus we have \eqref{eq:inf-sup-stability-Hu},
  but only when choosing $m=3$ the dual variable provides an error estimator.
  Also note, that we chose $S_{\underline H}^0(\Sigma)$ just for ease of
  presentation. It is sufficient to have a discrete space
  $\pwc_H^0(\Sigma)\subset X_{\underline H}\subset \pwc_h^0(\Sigma)$ that
  fulfills the discrete inf-sup stability \eqref{eq:inf-sup-stability-Hu}
  for all $v_{\underline H}\in X_{\underline H}$ and for which the solution
  $(\overline p_h,w_{\underline H})\in \pwc_h^0(\Sigma)\times X_{\underline H}$
  fulfills the saturation assumption \eqref{eq:saturation-assumption}.
\end{remark}

\section{Numerical Experiments}
\label{sec:NumExp}

\subsection{Set up}
We revisit two experiments, introduced in \cite{SteinbachUrzuaZank:2022}, and 
consider an additional experiment, all of which are posed on the same spatial 
domain $(0,3)$, i.e., $L=3$ and on the time interval $(0,6)$, i.e., $T=6$. 
We consider the following three different Dirichlet data
\begin{equation*}
    g_1(x,t):= \begin{cases}
      \frac{1}{2}(t-2)^3(-t)^3 & \mbox{for} \;
      0\leq t \leq 2\text{ and } x=0,\\
      \frac{1}{2}(t-5)^3(3-t)^3 & \mbox{for} \;
      L \leq t \leq L+2\text{ and } x=L,\\
        0 & \text{otherwise},
    \end{cases}
\end{equation*}
\begin{equation*}
    g_2(x,t):= \begin{cases}
      \frac{1}{2}|\sin(-\pi t)|& \mbox{for} \;
      0\leq t \text{ and } x=0,\\
      \frac{1}{2}|\sin(\pi(L-t))| & \mbox{for}
      \; L \leq t \text{ and } x=L,\\
        0 & \text{otherwise},
    \end{cases}
\end{equation*}
and
\begin{equation}
\label{def:g_c}
    g_3(x,t):= \begin{cases}
        t^{1/4} & \mbox{for} \; 0 \leq t \text{ and } x=0,\\
        (t-L)^{1/4} & \mbox{for} \; L \leq t \text{ and } x=L,\\
        0 & \text{otherwise}.
    \end{cases}
\end{equation}
We will be looking for solutions of the variational formulation for the mixed 
boundary element method as described in \eqref{mixed BEM}. For 
comparison, we will also consider two other variational formulations: 
the energetic BEM formulation, as described in
\cite{Aimi:2009, SteinbachUrzuaZank:2022}; 
and the \emph{modified Hilbert transform (MHT)} formulation from
\cite{SteinbachUrzuaZank:2022}. For 
clarity, let us restate the mixed boundary element methods for the specific 
numerical experiments, $i=1,2,3$. 
\begin{itemize}
\item \emph{standard least squares formulation:}
  find $p_h \in \pwc_h^0(\Sigma)\subset L^2(\Sigma)$ and
  $w_H \in \pwc_H^0(\Sigma)\subset [H^1_{,0}(\Sigma)]'$ such that
  \begin{equation}\label{vf:default}
    \langle p_h , q_h \rangle_{L^2(\Sigma)} +
    \langle \OV w_H , q_h \rangle_{L^2(\Sigma)} =
    \langle g_i , q_h \rangle_{L^2(\Sigma)}, \quad
    \langle p_h , \OV v_H \rangle_{L^2(\Sigma)} = 0,
  \end{equation}
  is satisfied for all $ q_h \in \pwc_h^0(\Sigma)$ and for all
  $v_H \in \pwc_H^0(\Sigma)$.
\item \emph{energetic BEM:}
  find $p_h \in \pwc_h^0(\Sigma)\subset L^2(\Sigma)$ and
  $z_H \in \pwc_H^0(\Sigma)\subset L^2(\Sigma)$ such that
  \begin{equation}\label{vf:EBEM}
    \langle p_h , q_h \rangle_{L^2(\Sigma)} +
    \langle \partial_t \OV z_H , q_h \rangle_{L^2(\Sigma)} =
    \langle \partial_t g_i , q_h \rangle_{L^2(\Sigma)}, \quad
    \langle p_h ,\partial_t \OV v_H \rangle_{L^2(\Sigma)} = 0,
  \end{equation}
  is satisfied for all $q_h \in \pwc_h^0(\Sigma)$ and for all
  $v_H \in \pwc_H^0(\Sigma)$.
\item \emph{modified Hilbert transform formulation:}
  Let $\mathcal{H}_T$ be the modified Hilbert transform defined 
  in \cite[eqn.~(2.8)]{SteinbachZank:2018}, 
  we want to solve: find $p_h \in \pwc_h^0(\Sigma)\subset L^2(\Sigma)$
  and $w_H \in \pwc_H^0(\Sigma) \subset [H_{,0}^1(\Sigma)]'$ such that
  \begin{equation}\label{vf:MHT}
    \langle p_h , q_h \rangle_{L^2(\Sigma)} +
    \langle \mathcal{H}_T \OV w_H , q_h \rangle_{L^2(\Sigma)} =
    \langle \mathcal{H}_T g_i , q_h \rangle_{L^2(\Sigma)}, \quad
    \langle p_h , \mathcal{H}_T \OV v_H \rangle_{L^2(\Sigma)} = 0,
  \end{equation}
  is satisfied for all $q_h \in \pwc_h^0(\Sigma)$ and for all
  $v_H \in \pwc_H^0(\Sigma)$. 
\end{itemize}
Throughout this section, the numerical experiments are implemented in Python.
For the solution of all linear systems built-in direct symmetric solvers
are used\footnote{\textsc{SciPy.Linalg.solve}}.

\subsection{Computation of the dual norm}\label{ssec: dual norm}
In order to compute the error $\|w-w_H\|_{[H^1_{,0}(\Sigma)]'}$, we require
the exact solution $w$, and a proper representation of the dual norm
$\|w-w_H\|_{[H^1_{,0}(\Sigma)]'}$.
In general, solutions to the indirect approach as considered in this
paper do not yield densities that can be interpreted physically. However,
in our specific setting, we are able to derive the exact density $w$ by
noting that for all functions $g_i$ we have
$g_i(0,t-L) = g_i(L,t)$ for $t\geq 0$.
We aim to find the exact solution $w_i$, satisfying $\OV w_i = g_i$.
Let us define 
\begin{equation*}
    \widetilde{w}_i(x,t):= \begin{cases}
        2 \partial_t\; g_i(0,t) & \mbox{for} \; x=0,\\
        0 & \mbox{for} \; x=L.
    \end{cases}
\end{equation*}
Then, one can verify that
\begin{equation*}
    \left(\OV  \widetilde{w}_i \right)(x,t) = \begin{cases}
        g_i(0,t) & \mbox{for} \; x= 0, \\
        g_i(0,t-L)=g_i(L,t) & \mbox{for} \; x=L,
    \end{cases}
  \end{equation*}
  i.e., $w_i=\widetilde{w}_i$.

In order to compute the dual norm $\|\cdot\|_{[H^1_{,0}(\Sigma)]'}$, note,
that by the Riesz representation theorem, for $w\in [H^1_{,0}(0,T)]'$
there exists exactly one $\phi_w\in H_{,0}^1(0,T)$ such that 
\begin{equation}\label{eq:norm-rep-VF}
  \langle \phi_w,v\rangle_{H^1_{,0}(0,T)} =
  \int_0^T \partial_t \phi_w(t)\partial_t v(t)\, dt =
  \langle w,v\rangle_{(0,T)}\quad  \text{for all }v\in H^1_{,0}(0,T),
\end{equation}
and that 
\begin{equation}
  \norm{\phi_w}_{H^1_{,0}(0,T)}^2 =
  \norm{\partial_t \phi_w}_{L^2(0,T)}^2 =
  \langle w,\phi_w\rangle_{(0,T)} =\norm{w}_{[H_{,0}^1(0,T)]'}^2.
\end{equation}
Note that \eqref{eq:norm-rep-VF} is the variational formulation of the
boundary value problem
\begin{equation*}
  -\partial_{tt}\phi_w(t) = w(t) \mbox{ for } t\in(0,T),
  \quad \partial_t\phi_w(0)=0,\quad \phi_w(T) =0, 
\end{equation*}
for which the solution is given, using Greens function, as
\begin{equation*}
  \phi_w(t) = \int_0^T G(t,s)w(s)\, ds, \quad
  \text{where } G(t,s) =
  \begin{cases}
    T-t,& s\in (0,t),\\
    T-s,& s\in (t,T).
  \end{cases}
\end{equation*} 

\subsection{Numerical results}
We start by checking numerically if the theoretical results in
Theorems \ref{Thm Sh} and \ref{thm:discrete-stability-n-time-slices},
and Lemma~\ref{lem:aPriori} are sharp in excluding
$m=2$ (namely, when each element of $\Sigma_H$ is decomposed into
two equally sized elements to obtain the fine mesh $\Sigma_h$). 
Fig.~\ref{Fig:natural_m2_AB} shows the results for the standard
formulation \eqref{vf:default}, given $m=2$, while those of 
energetic BEM and MHT are displayed in Fig.~\ref{Fig:EBEM_m2_AB} and
\ref{Fig:MHT_m2_AB}, respectively. In all cases, the method converges.
However, we see that they behave differently when considering adaptive
refinements.

It is clear from Fig.~\ref{Fig:natural_m2_AB} that $\norm{p_h}_{L^2(\Sigma)}$
does not provide a reliable error estimator for the standard formulation
\eqref{vf:default} when $m=2$. This fits the theory presented in
Lemma~\ref{Lem ph error indicator}, which states that, in order to show
$\norm{p_h}_{L^2(\Sigma)}$ is an error estimator, the saturation assumption
\eqref{eq:saturation-assumption} must hold. This only happens when $m>2$
for the standard formulation \eqref{vf:default}, while it is already true
for $m=2$ for energetic BEM and MHT, since these formulations are discrete
inf-sup stable for the case $h=H$.

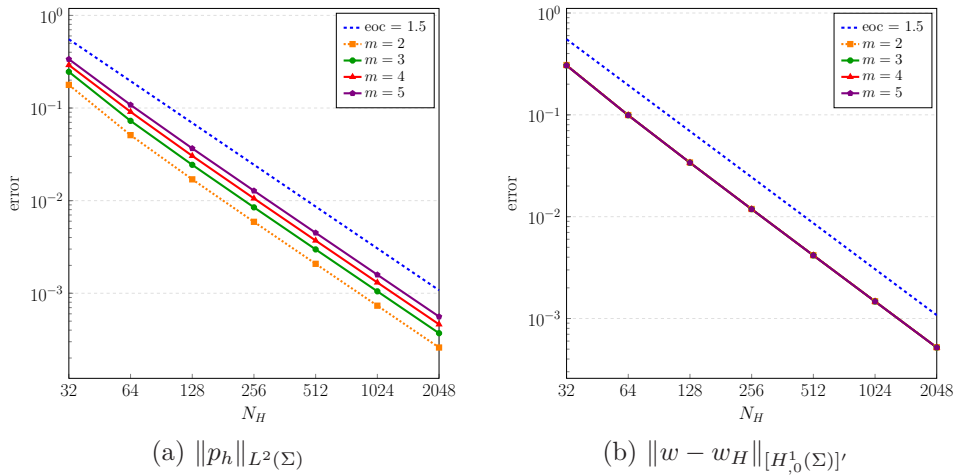
\begin{figure}[h]
    \centering
    \subfloat[][\centering  $\|p_h\|_{L^2(\Sigma)}$]{{\pgfplotstabletranspose\two{Data/default_uni_B_m2.txt}
\pgfplotstabletranspose\three{Data/default_uni_B_m3.txt}
\pgfplotstabletranspose\four{Data/default_uni_B_m4.txt}
\pgfplotstabletranspose\five{Data/default_uni_B_m5.txt}


\pgfplotstableset{
    create on use/Nh/.style={
        create col/copy column from table={Data/N.txt}{N}}
    }
{\scalefont{1.5}
\begin{tikzpicture}
\begin{loglogaxis}[
    width=14cm,
    height=14cm,
    xlabel={$N_H$},
    ylabel={error},
    xmin=32, xmax=2048,
    ytick pos= left,
    xtick pos= bottom,
    xtick={32,64,128,256,512,1024,2048},
    xticklabels= {32,64,128,256,512,1024,2048},
    legend pos=north east,
    legend cell align={left},
    legend style ={font = \large, inner sep= 6pt, outer sep = -6pt},
    ymajorgrids=true,
    grid style={dashed, line width=0.1pt, lightgray!50},
]

\addplot [
    domain= 32:2048, 
    samples=8, 
    color=blue,
    dashed,
    line width=2pt
]
{100*x^-(1.5)};
\addlegendentry{eoc = 1.5}


\addplot[orange,mark= square*,line width=2pt, mark options={solid},dotted] table [x = Nh, y index = 1] {\two};
\addlegendentry{$m=2$}

\addplot[green!60!black, mark=oplus*,mark options=solid, line width=2pt] table [x = Nh, y index = 1] {\three};
\addlegendentry{$m=3$}

\addplot[red, mark=triangle*,mark options=solid, line width=2pt] table [x = Nh, y index = 1] {\four};
\addlegendentry{$m=4$}

\addplot[violet, mark=pentagon*,mark options=solid, line width=2pt] table [x = Nh, y index = 1] {\five};
\addlegendentry{$m=5$}

\end{loglogaxis}
\end{tikzpicture}
} }}%
    \subfloat[][\centering $\|w-w_H\|_{[H^1_{,0}(\Sigma)]'}$]{{\pgfplotstabletranspose\two{Data/default_uni_B_m2.txt}
\pgfplotstabletranspose\three{Data/default_uni_B_m3.txt}
\pgfplotstabletranspose\four{Data/default_uni_B_m4.txt}
\pgfplotstabletranspose\five{Data/default_uni_B_m5.txt}


\pgfplotstableset{
    create on use/Nh/.style={
        create col/copy column from table={Data/N.txt}{N}}
    }
{\scalefont{1.5}
\begin{tikzpicture}
\begin{loglogaxis}[
    width=14cm,
    height=14cm,
    xlabel={$N_H$},
    ylabel={error},
    xmin=32, xmax=2048,
    ytick pos= left,
    xtick pos= bottom,
    xtick={32,64,128,256,512,1024,2048},
    xticklabels= {32,64,128,256,512,1024,2048},
    legend pos=north east,
    legend cell align={left},
    legend style ={font = \large, inner sep= 6pt, outer sep = -6pt},
    ymajorgrids=true,
    grid style={dashed, line width=0.1pt, lightgray!50},
]

\addplot [
    domain= 32:2048, 
    samples=8, 
    color=blue,
    dashed,
    line width=2pt
]
{100*x^-(1.5)};
\addlegendentry{eoc = 1.5}


\addplot[orange,mark= square*,line width=2pt, mark options={solid}, dotted] table [x = Nh, y index = 2] {\two};
\addlegendentry{$m=2$}

\addplot[green!60!black, mark=oplus*,mark options=solid, line width=2pt] table [x = Nh, y index = 2] {\three};
\addlegendentry{$m=3$}

\addplot[red, mark=triangle*,mark options=solid, line width=2pt] table [x = Nh, y index = 2] {\four};
\addlegendentry{$m=4$}

\addplot[violet, mark=pentagon*,mark options=solid, line width=2pt] table [x = Nh, y index = 2] {\five};
\addlegendentry{$m=5$}

\end{loglogaxis}
\end{tikzpicture}
} }}%
    \caption{The error indicator and the error of the solution when using the standard formulation \eqref{vf:default}, uniform meshes,
    Dirichlet datum $g_2$, and different choices of $m$.}%
    \label{Fig:error_uniform_1a1b_all}%
\end{figure}

In order to verify that $\norm{p_h}_{L^2(\Sigma)}$ becomes an error estimator
for the standard formulation \eqref{vf:default} when $m>2$,
Fig.~\ref{Fig:natural_m3_AB} depicts the results for the standard
formulation \eqref{vf:default} when we consider a fine mesh $\Sigma_h$ such
that each element of $\Sigma_H$ is decomposed into three $(m=3)$
equally sized elements. Interestingly, in the uniform case the exact
error $\|w-w_H\|_{[H^1_{,0}(\Sigma)]'}$ does not change significantly,
\textit{but} the convergence rate of $\|p_h\|_{L^2(\Sigma)}$ seems to
correspond to the convergence rate of the exact error in this case. 
As shown in Fig.~\ref{Fig:error_uniform_1a1b_all}, further increasing the
value of $m$ does not seem to affect the convergence rate of the the error
indicator $\|p_h\|_{L^2(\Sigma)}$ and the error $\|w-w_H\|_{[H^1_{,0}(\Sigma)]'}$
for uniform refinements of the standard formulation \eqref{vf:default}.
Moreover, the difference in the error $\|w-w_H\|_{[H^1_{,0}(\Sigma)]'}$ seems to 
be negligible for different choices of $m$.

\begin{figure}[h]
    \centering
    \subfloat[][\centering $g_1$]{{\pgfplotstabletranspose\uni{Data/default_uni_A_m2.txt}
\pgfplotstabletranspose\adapt{Data/default_adapt_A_m2.txt}


\pgfplotstableset{
    create on use/Nh/.style={
        create col/copy column from table={Data/N.txt}{N}}
    }

{\scalefont{1.5}
\begin{tikzpicture}
\begin{loglogaxis}[
    width=14cm,
    height=14cm,
    xlabel={$N_H$},
    ylabel={error},
    xmin=32, xmax=2048,
    ytick pos= left,
    xtick pos= bottom,
    xtick={32,64,128,256,512,1024,2048},
    xticklabels= {32,64,128,256,512,1024,2048},
    legend pos=north east,
    legend cell align={left},
    legend style ={font = \large, inner sep= 6pt, outer sep = -6pt},
    ymajorgrids=true,
    grid style={dashed, line width=0.1pt, lightgray!50},
]

\addplot [
    domain= 32:2048, 
    samples=8, 
    color=blue,
    dashed,
    line width=2pt
]
{100*x^-2};
\addlegendentry{eoc = 2}


\addplot[orange,mark= square*,line width=2pt, mark options={solid}, dotted] table [x = Nh, y index = 1] {\uni};
\addlegendentry{$\|p_h\|_{L^2(\Sigma)}$-uni}

\addplot[green!60!black, mark=square*,mark options=solid, line width=2pt] table [x = Nh, y index = 2] {\uni};
\addlegendentry{$\|w-w_H\|_{[H_{,0}^1(\Sigma)]^\prime}$-uni}

\addplot[red,mark= triangle*,line width=2pt, mark options={solid}, dotted] table [x index=3, y index = 1] {\adapt};
\addlegendentry{$\|p_h\|_{L^2(\Sigma)}$-ad}

\addplot[violet, mark=triangle*,mark options=solid, line width=2pt] table [x index= 3, y index = 2] {\adapt};
\addlegendentry{$\|w-w_H\|_{[H_{,0}^1(\Sigma)]^\prime}$-ad}

\end{loglogaxis}
\end{tikzpicture}
} }}%
    \subfloat[][\centering $g_2$]{{\pgfplotstabletranspose\uni{Data/default_uni_B_m2.txt}
\pgfplotstabletranspose\adapt{Data/default_adapt_B_m2.txt}


\pgfplotstableset{
    create on use/Nh/.style={
        create col/copy column from table={Data/N.txt}{N}}
    }

{\scalefont{1.5}
\begin{tikzpicture}
\begin{loglogaxis}[
    width=14cm,
    height=14cm,
    xlabel={$N_H$},
    ylabel={error},
    xmin=32, xmax=2048,
    ytick pos= left,
    xtick pos= bottom,
    xtick={32,64,128,256,512,1024,2048},
    xticklabels= {32,64,128,256,512,1024,2048},
    legend pos=north east,
    legend cell align={left},
    legend style ={font = \large, inner sep= 6pt, outer sep = -6pt},
    ymajorgrids=true,
    grid style={dashed, line width=0.1pt, lightgray!50},
]

\addplot [
    domain= 32:2048, 
    samples=8, 
    color=blue,
    dashed,
    line width=2pt
]
{100*x^-(1.5)};
\addlegendentry{eoc = 1.5}


\addplot[orange,mark= square*,line width=2pt, mark options={solid}, dotted] table [x = Nh, y index = 1] {\uni};
\addlegendentry{$\|p_h\|_{L^2(\Sigma)}$-uni}

\addplot[green!60!black, mark=square*,mark options=solid, line width=2pt] table [x = Nh, y index = 2] {\uni};
\addlegendentry{$\|w-w_H\|_{[H_{,0}^1(\Sigma)]^\prime}$-uni}

\addplot[red,mark= triangle*,line width=2pt, mark options={solid}, dotted] table [x index=3, y index = 1] {\adapt};
\addlegendentry{$\|p_h\|_{L^2(\Sigma)}$-ad}

\addplot[violet, mark=triangle*,mark options=solid, line width=2pt] table [x index= 3, y index = 2] {\adapt};
\addlegendentry{$\|w-w_H\|_{[H_{,0}^1(\Sigma)]^\prime}$-ad}

\end{loglogaxis}
\end{tikzpicture}
} }}%
    \caption{
    Comparison of errors and error indicators 
    for uniform and adaptive refinement using the standard formulation \eqref{vf:default}, $m = 2$,
    and different Dirichlet data.}%
    \label{Fig:natural_m2_AB}%
\end{figure}
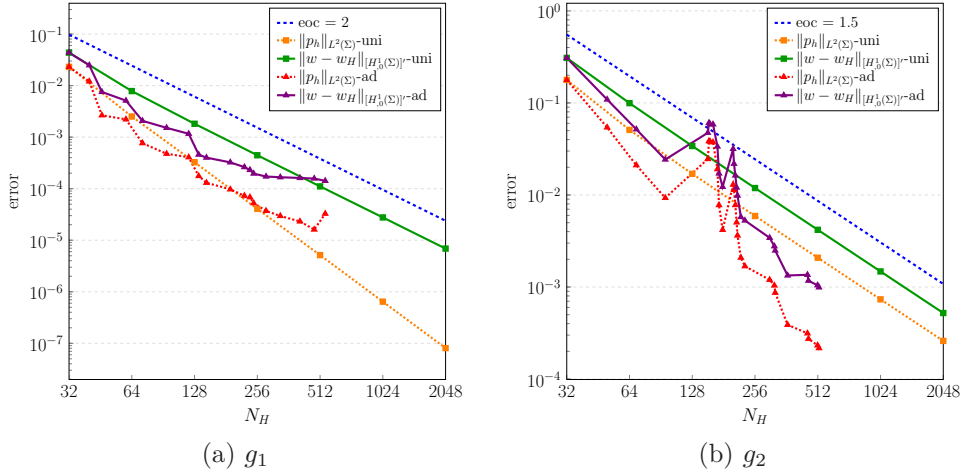
\begin{figure}[h]
    \centering
    \subfloat[][\centering $g_1$]{{\pgfplotstabletranspose\uni{Data/EBEM_uni_A_m2.txt}
\pgfplotstabletranspose\adapt{Data/EBEM_adapt_A_m2.txt}


\pgfplotstableset{
    create on use/Nh/.style={
        create col/copy column from table={Data/N.txt}{N}}
    }
{\scalefont{1.5}
\begin{tikzpicture}
\begin{loglogaxis}[
    width=14cm,
    height=14cm,
    xlabel={$N_H$},
    ylabel={error},
    xmin=32, xmax=2048,
    ytick pos= left,
    xtick pos= bottom,
    xtick={32,64,128,256,512,1024,2048},
    xticklabels= {32,64,128,256,512,1024,2048},
    legend pos=north east,
    legend cell align={left},
    legend style ={font = \large, inner sep= 6pt, outer sep = -6pt},
    ymajorgrids=true,
    grid style={dashed, line width=0.1pt, lightgray!50},
]

\addplot [
    domain= 32:2048, 
    samples=8, 
    color=blue,
    dashed,
    line width=2pt
]
{30*x^-1};
\addlegendentry{eoc = 1}


\addplot[orange,mark= square*,line width=2pt, mark options={solid}, dotted] table [x = Nh, y index = 1] {\uni};
\addlegendentry{$\|p_h\|_{L^2(\Sigma)}$-uni}

\addplot[green!60!black, mark=square*,mark options=solid, line width=2pt] table [x = Nh, y index = 2] {\uni};
\addlegendentry{$\|w-w_H\|_{L^2(\Sigma)}$-uni}

\addplot[red,mark= triangle*,line width=2pt, mark options={solid}, dotted] table [x index=3, y index = 1] {\adapt};
\addlegendentry{$\|p_h\|_{L^2(\Sigma)}$-ad}

\addplot[violet, mark=triangle*,mark options=solid, line width=2pt] table [x index= 3, y index = 2] {\adapt};
\addlegendentry{$\|w-w_H\|_{L^2(\Sigma)}$-ad}

\end{loglogaxis}
\end{tikzpicture}
} }}%
    \subfloat[][\centering $g_2$]{{\pgfplotstabletranspose\uni{Data/EBEM_uni_B_m2.txt}
\pgfplotstabletranspose\adapt{Data/EBEM_adapt_B_m2.txt}


\pgfplotstableset{
    create on use/Nh/.style={
        create col/copy column from table={Data/N.txt}{N}}
    }
{\scalefont{1.5}
\begin{tikzpicture}
\begin{loglogaxis}[
    width=14cm,
    height=14cm,
    xlabel={$N_H$},
    ylabel={error},
    xmin=32, xmax=2048,
    ytick pos= left,
    xtick pos= bottom,
    xtick={32,64,128,256,512,1024,2048},
    xticklabels= {32,64,128,256,512,1024,2048},
    legend pos=north east,
    legend cell align={left},
    legend style ={font = \large, inner sep= 6pt, outer sep = -6pt},
    ymajorgrids=true,
    grid style={dashed, line width=0.1pt, lightgray!50},
]

\addplot [
    domain= 32:2048, 
    samples=8, 
    color=blue,
    dashed,
    line width=2pt
]
{30*x^-0.5};
\addlegendentry{eoc = 0.5}


\addplot[orange,mark= square*,line width=2pt, mark options={solid}, dotted] table [x = Nh, y index = 1] {\uni};
\addlegendentry{$\|p_h\|_{L^2(\Sigma)}$-uni}

\addplot[green!60!black, mark=square*,mark options=solid, line width=2pt] table [x = Nh, y index = 2] {\uni};
\addlegendentry{$\|w-w_H\|_{L^2(\Sigma)}$-uni}

\addplot[red,mark= triangle*,line width=2pt, mark options={solid}, dotted] table [x index=3, y index = 1] {\adapt};
\addlegendentry{$\|p_h\|_{L^2(\Sigma)}$-ad}

\addplot[violet, mark=triangle*,mark options=solid, line width=2pt] table [x index= 3, y index = 2] {\adapt};
\addlegendentry{$\|w-w_H\|_{L^2(\Sigma)}$-ad}

\end{loglogaxis}
\end{tikzpicture}
} }}%
    \caption{
    Comparison of errors and error indicators 
    for uniform and adaptive refinement using the energetic formulation \eqref{vf:EBEM}, $m = 2$, 
    and different Dirichlet data.}%
    \label{Fig:EBEM_m2_AB}%
\end{figure}
\begin{figure}[h]
    \centering
    \subfloat[][\centering $g_1$]{{\pgfplotstabletranspose\uni{Data/MHT_uni_A_m2.txt}
\pgfplotstabletranspose\adapt{Data/MHT_adapt_A_m2.txt}


\pgfplotstableset{
    create on use/Nh/.style={
        create col/copy column from table={Data/N.txt}{N}}
    }
{\scalefont{1.5}
\begin{tikzpicture}
\begin{loglogaxis}[
    width=14cm,
    height=14cm,
    xlabel={$N_H$},
    ylabel={error},
    xmin=32, xmax=2048,
    ytick pos= left,
    xtick pos= bottom,
    xtick={32,64,128,256,512,1024,2048},
    xticklabels= {32,64,128,256,512,1024,2048},
    legend pos=north east,
    legend cell align={left},
    legend style ={font = \large, inner sep= 6pt, outer sep = -6pt},
    ymajorgrids=true,
    grid style={dashed, line width=0.1pt, lightgray!50},
]

\addplot [
    domain= 32:2048, 
    samples=8, 
    color=blue,
    dashed,
    line width=2pt
]
{100*x^-2};
\addlegendentry{eoc = 2}


\addplot[orange,mark= square*,line width=2pt, mark options={solid}, dotted] table [x = Nh, y index = 1] {\uni};
\addlegendentry{$\|p_h\|_{L^2(\Sigma)}$-uni}

\addplot[green!60!black, mark=square*,mark options=solid, line width=2pt] table [x = Nh, y index = 2] {\uni};
\addlegendentry{$\|w-w_H\|_{[H_{,0}^1(\Sigma)]^\prime}$-uni}

\addplot[red,mark= triangle*,line width=2pt, mark options={solid}, dotted] table [x index=3, y index = 1] {\adapt};
\addlegendentry{$\|p_h\|_{L^2(\Sigma)}$-ad}

\addplot[violet, mark=triangle*,mark options=solid, line width=2pt] table [x index= 3, y index = 2] {\adapt};
\addlegendentry{$\|w-w_H\|_{[H_{,0}^1(\Sigma)]^\prime}$-ad}

\end{loglogaxis}
\end{tikzpicture}
} }}%
    \subfloat[][\centering $g_2$]{{\pgfplotstabletranspose\uni{Data/MHT_uni_B_m2.txt}
\pgfplotstabletranspose\adapt{Data/MHT_adapt_B_m2.txt}


\pgfplotstableset{
    create on use/Nh/.style={
        create col/copy column from table={Data/N.txt}{N}}
    }

{\scalefont{1.5}
\begin{tikzpicture}
\begin{loglogaxis}[
    width=14cm,
    height=14cm,
    xlabel={$N_H$},
    ylabel={error},
    xmin=32, xmax=2048,
    ytick pos= left,
    xtick pos= bottom,
    xtick={32,64,128,256,512,1024,2048},
    xticklabels= {32,64,128,256,512,1024,2048},
    legend pos=north east,
    legend cell align={left},
    legend style ={font = \large, inner sep= 6pt, outer sep = -6pt},
    ymajorgrids=true,
    grid style={dashed, line width=0.1pt, lightgray!50},
]

\addplot [
    domain= 32:2048, 
    samples=8, 
    color=blue,
    dashed,
    line width=2pt
]
{100*x^-1.5};
\addlegendentry{eoc = 1.5}


\addplot[orange,mark= square*,line width=2pt, mark options={solid}, dotted] table [x = Nh, y index = 1] {\uni};
\addlegendentry{$\|p_h\|_{L^2(\Sigma)}$-uni}

\addplot[green!60!black, mark=square*,mark options=solid, line width=2pt] table [x = Nh, y index = 2] {\uni};
\addlegendentry{$\|w-w_H\|_{[H_{,0}^1(\Sigma)]^\prime}$-uni}

\addplot[red,mark= triangle*,line width=2pt, mark options={solid}, dotted] table [x index=3, y index = 1] {\adapt};
\addlegendentry{$\|p_h\|_{L^2(\Sigma)}$-ad}

\addplot[violet, mark=triangle*,mark options=solid, line width=2pt] table [x index= 3, y index = 2] {\adapt};
\addlegendentry{$\|w-w_H\|_{[H_{,0}^1(\Sigma)]^\prime}$-ad}

\end{loglogaxis}
\end{tikzpicture}
} }}%
    \caption{
    Comparison of errors and error indicators 
    for uniform and adaptive refinement using the MHT formulation \eqref{vf:MHT}, $m = 2$, 
    and different Dirichlet data. }%
    \label{Fig:MHT_m2_AB}%
\end{figure}
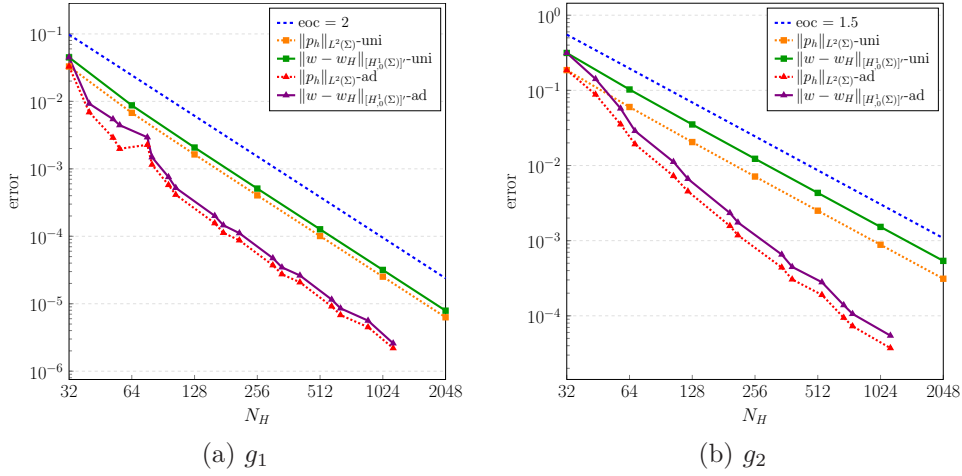

\begin{figure}[h]
    \centering
    \subfloat[][\centering $g_1$]{{\pgfplotstabletranspose\uni{Data/default_uni_A_m3.txt}
\pgfplotstabletranspose\adapt{Data/default_adapt_A_m3.txt}


\pgfplotstableset{
    create on use/Nh/.style={
        create col/copy column from table={Data/N.txt}{N}}
    }
{\scalefont{1.5}
\begin{tikzpicture}
\begin{loglogaxis}[
    width=14cm,
    height=14cm,
    xlabel={$N_H$},
    ylabel={error},
    xmin=32, xmax=2048,
    ytick pos= left,
    xtick pos= bottom,
    xtick={32,64,128,256,512,1024,2048},
    xticklabels= {32,64,128,256,512,1024,2048},
    legend pos=north east,
    legend cell align={left},
    legend style ={font = \large, inner sep= 6pt, outer sep = -6pt},
    ymajorgrids=true,
    grid style=dashed,
]

\addplot [
    domain= 32:2048, 
    samples=8, 
    color=blue,
    dashed,
    line width=2pt
]
{100*x^-2};
\addlegendentry{eoc = 2}


\addplot[orange,mark= square*,line width=2pt, mark options={solid}, dotted] table [x = Nh, y index = 1] {\uni};
\addlegendentry{$\|p_h\|_{L^2(\Sigma)}$-uni}

\addplot[green!60!black, mark=square*,mark options=solid, line width=2pt] table [x = Nh, y index = 2] {\uni};
\addlegendentry{$\|w-w_H\|_{[H_{,0}^1(\Sigma)]^\prime}$-uni}

\addplot[red,mark= triangle*,line width=2pt, mark options={solid}, dotted] table [x index=3, y index = 1] {\adapt};
\addlegendentry{$\|p_h\|_{L^2(\Sigma)}$-ad}

\addplot[violet, mark=triangle*,mark options=solid, line width=2pt] table [x index= 3, y index = 2] {\adapt};
\addlegendentry{$\|w-w_H\|_{[H_{,0}^1(\Sigma)]^\prime}$-ad}

\end{loglogaxis}
\end{tikzpicture}
} }}%
    \subfloat[][\centering $g_2$]{{\pgfplotstabletranspose\uni{Data/default_uni_B_m3.txt}
\pgfplotstabletranspose\adapt{Data/default_adapt_B_m3.txt}


\pgfplotstableset{
    create on use/Nh/.style={
        create col/copy column from table={Data/N.txt}{N}}
    }
{\scalefont{1.5}
\begin{tikzpicture}
\begin{loglogaxis}[
    width=14cm,
    height=14cm,
    xlabel={$N_H$},
    ylabel={error},
    xmin=32, xmax=2048,
    ytick pos= left,
    xtick pos= bottom,
    xtick={32,64,128,256,512,1024,2048},
    xticklabels= {32,64,128,256,512,1024,2048},
    legend pos=north east,
    legend cell align={left},
    legend style ={font = \large, inner sep= 6pt, outer sep = -6pt},
    ymajorgrids=true,
    grid style={dashed, line width=0.1pt, lightgray!50},
]

\addplot [
    domain= 32:2048, 
    samples=8, 
    color=blue,
    dashed,
    line width=2pt,
]
{100*x^-(1.5)};
\addlegendentry{eoc = 1.5}


\addplot[orange,mark= square*,line width=2pt, mark options={solid}, dotted] table [x = Nh, y index = 1] {\uni};
\addlegendentry{$\|p_h\|_{L^2(\Sigma)}$-uni}

\addplot[green!60!black, mark=square*,mark options=solid, line width=2pt] table [x = Nh, y index = 2] {\uni};
\addlegendentry{$\|w-w_H\|_{[H_{,0}^1(\Sigma)]^\prime}$-uni}

\addplot[red,mark= triangle*,line width=2pt, mark options={solid}, dotted] table [x index=3, y index = 1] {\adapt};
\addlegendentry{$\|p_h\|_{L^2(\Sigma)}$-ad}

\addplot[violet, mark=triangle*,mark options=solid, line width=2pt] table [x index= 3, y index = 2] {\adapt};
\addlegendentry{$\|w-w_H\|_{[H_{,0}^1(\Sigma)]^\prime}$-ad}

\end{loglogaxis}
\end{tikzpicture}
} }}%
    \caption{Comparison of errors and error indicators 
    for uniform and adaptive refinement using the standard formulation \eqref{vf:default}, $m = 3$,
    and different Dirichlet data.}%
    \label{Fig:natural_m3_AB}%
\end{figure}

\begin{figure}[h]%
 \centering
 \subfloat[Standard formulation \eqref{vf:default}]{\pgfplotstabletranspose\uni{Data/default_uni_C_m2.txt}
\pgfplotstabletranspose\adapt{Data/default_adapt_C_m3.txt}


\pgfplotstableset{
    create on use/Nh/.style={
        create col/copy column from table={Data/N.txt}{N}}
    }
{\scalefont{1.5}
\begin{tikzpicture}
\begin{loglogaxis}[
    width=14cm,
    height=14cm,
    xlabel={$N_H$},
    ylabel={error},
    xmin=32, xmax=2048,
    ytick pos= left,
    xtick pos= bottom,
    xtick={32,64,128,256,512,1024,2048},
    xticklabels= {32,64,128,256,512,1024,2048},
    legend pos=north east,
    legend cell align={left},
    legend style ={font = \large, inner sep= 6pt, outer sep = -6pt},
    ymajorgrids=true,
    grid style={dashed, line width=0.1pt, lightgray!50},
]

\addplot [
    domain= 32:2048, 
    samples=8, 
    color=blue,
    dashed,
    line width=2pt,
]
{30*x^-(0.75)};
\addlegendentry{eoc =0.75}


\addplot[orange,mark= square*,line width=2pt, mark options={solid}, dotted] table [x = Nh, y index = 1] {\uni};
\addlegendentry{$\|p_h\|_{L^2(\Sigma)}$-uni}

\addplot[green!60!black, mark=square*,mark options=solid, line width=2pt] table [x = Nh, y index = 2] {\uni};
\addlegendentry{$\|w-w_H\|_{[H_{,0}^1(\Sigma)]^\prime}$-uni}

\addplot[red,mark= triangle*,line width=2pt, mark options={solid}, dotted] table [x index=3, y index = 1] {\adapt};
\addlegendentry{$\|p_h\|_{L^2(\Sigma)}$-ad}

\addplot[violet, mark=triangle*,mark options=solid, line width=2pt] table [x index= 3, y index = 2] {\adapt};
\addlegendentry{$\|w-w_H\|_{[H_{,0}^1(\Sigma)]^\prime}$-ad}

\end{loglogaxis}
\end{tikzpicture}
} }%
 \subfloat[Energetic BEM \eqref{vf:EBEM}]{\pgfplotstabletranspose\uni{Data/EBEM_uni_C_m2.txt}
\pgfplotstabletranspose\adapt{Data/EBEM_adapt_C_m2.txt}


\pgfplotstableset{
    create on use/Nh/.style={
        create col/copy column from table={Data/N.txt}{N}}
    }

{\scalefont{1.5}
\begin{tikzpicture}
\begin{loglogaxis}[
    width=14cm,
    height=14cm,
    xlabel={$N_H$},
    ylabel={error},
    xmin=32, xmax=2048,
    ytick pos= left,
    xtick pos= bottom,
    xtick={32,64,128,256,512,1024,2048},
    xticklabels= {32,64,128,256,512,1024,2048},
    legend pos=north east,
    legend cell align={left},
    legend style ={font = \large, inner sep= 6pt, outer sep = -6pt},
    ymajorgrids=true,
    grid style={dashed, line width=0.1pt, lightgray!50},
]

\addplot [
    domain= 32:2048, 
    samples=8, 
    color=blue,
    dashed,
    line width=2pt
]
{30*x^-0.75};
\addlegendentry{eoc = 0.75}


\addplot[orange,mark= square*,line width=2pt, mark options={solid}, dotted] table [x = Nh, y index = 1] {\uni};
\addlegendentry{$\|p_h\|_{L^2(\Sigma)}$-uni}

\addplot[green!60!black, mark=square*,mark options=solid, line width=2pt] table [x = Nh, y index = 2] {\uni};
\addlegendentry{$\|w-w_H\|_{[H_{,0}^1(\Sigma)]^\prime}$-uni}

\addplot[red,mark= triangle*,line width=2pt, mark options={solid}, dotted] table [x index=3, y index = 1] {\adapt};
\addlegendentry{$\|p_h\|_{L^2(\Sigma)}$-ad}

\addplot[violet, mark=triangle*,mark options=solid, line width=2pt] table [x index= 3, y index = 2] {\adapt};
\addlegendentry{$\|w-w_H\|_{[H_{,0}^1(\Sigma)]^\prime}$-ad}

\end{loglogaxis}
\end{tikzpicture}
}}\\
 \subfloat[MHT formulation \eqref{vf:MHT}]{\pgfplotstabletranspose\uni{Data/MHT_uni_C_m2.txt}
\pgfplotstabletranspose\adapt{Data/MHT_adapt_C_m2.txt}


\pgfplotstableset{
    create on use/Nh/.style={
        create col/copy column from table={Data/N.txt}{N}}
    }

{\scalefont{1.5}
\begin{tikzpicture}
\begin{loglogaxis}[
    width=14cm,
    height=14cm,
    xlabel={$N_H$},
    ylabel={error},
    xmin=32, xmax=2048,
    ytick pos= left,
    xtick pos= bottom,
    xtick={32,64,128,256,512,1024,2048},
    xticklabels= {32,64,128,256,512,1024,2048},
    legend pos=north east,
    legend cell align={left},
    legend style ={font = \large, inner sep= 6pt, outer sep = -6pt},
    ymajorgrids=true,
    grid style={dashed, line width=0.1pt, lightgray!50},
]

\addplot [
    domain= 32:2048, 
    samples=8, 
    color=blue,
    dashed,
    line width=2pt
]
{30*x^-0.75};
\addlegendentry{eoc = 0.75}


\addplot[orange,mark= square*,line width=2pt, mark options={solid}, dotted] table [x = Nh, y index = 1] {\uni};
\addlegendentry{$\|p_h\|_{L^2(\Sigma)}$-uni}

\addplot[green!60!black, mark=square*,mark options=solid, line width=2pt] table [x = Nh, y index = 2] {\uni};
\addlegendentry{$\|w-w_H\|_{[H_{,0}^1(\Sigma)]^\prime}$-uni}

\addplot[red,mark= triangle*,line width=2pt, mark options={solid}, dotted] table [x index=3, y index = 1] {\adapt};
\addlegendentry{$\|p_h\|_{L^2(\Sigma)}$-ad}

\addplot[violet, mark=triangle*,mark options=solid, line width=2pt] table [x index= 3, y index = 2] {\adapt};
\addlegendentry{$\|w-w_H\|_{[H_{,0}^1(\Sigma)]^\prime}$-ad}

\end{loglogaxis}
\end{tikzpicture}
}}%
 \caption{
 Comparison of errors and error indicators for uniform and adaptive refinement schemes given Dirichlet datum $g_3$.}%
 \label{fig:all_g_c}%
\end{figure}
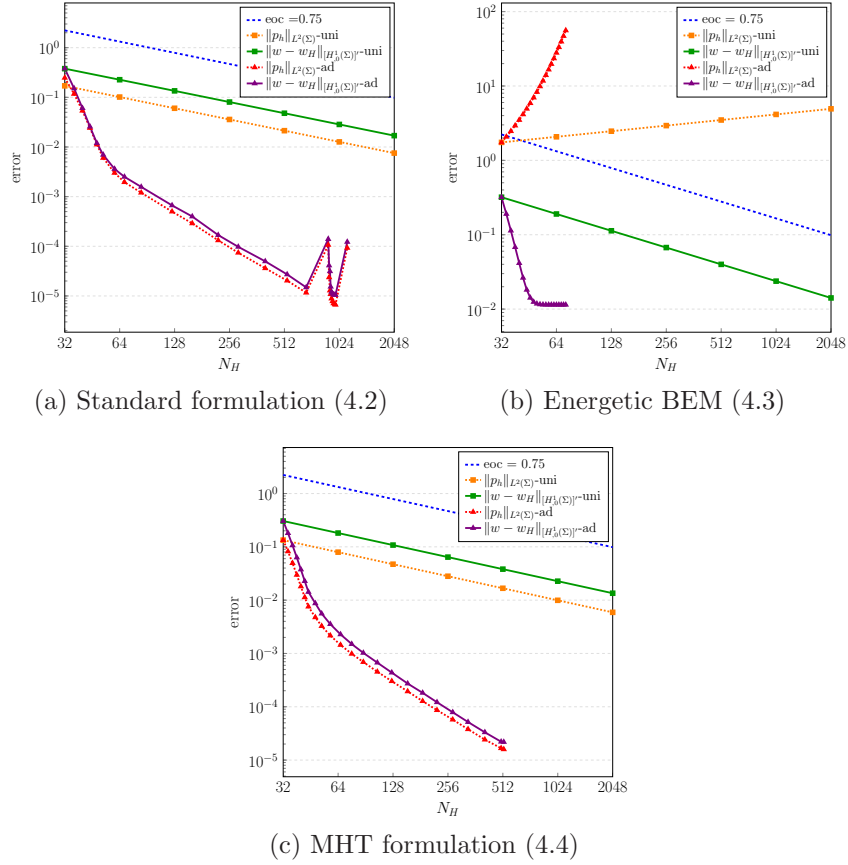

Up until this point, we have not yet considered results related to $g_3$.
For notational convenience, let us define $w$ as the exact solution to
the BIE for either the standard, energetic or MHT formulation. Then, given
Dirichlet data $g_3$ as defined in \eqref{def:g_c}, the density $w$ will
\emph{not} be in $L^2(\Sigma)$. Hence, solutions of this kind do not fit
our current framework for the energetic formulation \eqref{vf:EBEM}. With
that being said, we still adhere to the same energetic formulation and
discretization as already considered for $g_1$ and $g_2$. When it comes to
energetic formulations with $g_3$, only the norm in which the error is
measured is changed into $[H^1_{,0}(\Sigma)]^\prime$ as opposed to the
usual $L^2(\Sigma)$ norm. The results related to $g_3$, for different
formulations, are presented in Fig.~\ref{fig:all_g_c}. There we see that the
three different formulations converge with rate $0.75$ on uniform meshes.
For energetic BEM, $\norm{p_h}_{L^2(\Sigma)}$ no longer serves as an error
estimator for the error in the $[H^1_{,0}(\Sigma)]^\prime$ norm. This
explains why convergence of the adaptive routine for energetic BEM halts
after some refinements.

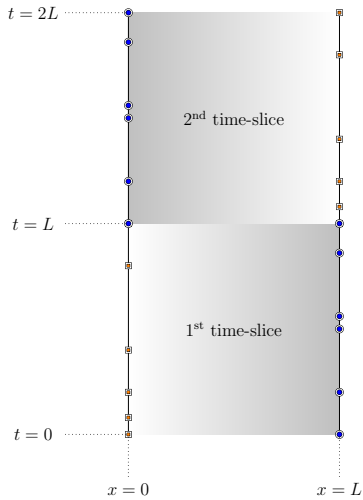
\begin{figure}[H]
\vspace{-.25cm}
    \centering
    \begin{tikzpicture}

\shade[left color=white,right color=lightgray] (0,0) rectangle (5,5) node[pos=.5] {$1^{\text{st}}$ time-slice};
\shade[left color=lightgray,right color=white] (0,5) rectangle (5,10) node[pos=.5] {$2^{\text{nd}}$ time-slice};
\def\lastx{-0.05}
\draw[black] (0,0)--(0,10);
\draw[black] (5,0)--(5,10);
\foreach \x[remember =\x as \lastx] in {0,0.4, 1, 2, 4 ,5}{
    \draw (0,\x) node[rectangle, fill =orange, scale = 0.4, draw=black, double] {};
    \draw (5,\x+5) node[rectangle, fill =orange, scale = 0.4, draw= black,double] {};
    }
\def\lastx{-0.05}
\foreach \x[remember =\x as \lastx] in {0,1,2.5,2.8,4.3,5}{
    \draw (0,\x+5) node[circle, fill =blue, scale = 0.4, draw= black, double] {};
    \draw (5,\x) node[circle, fill =blue, scale = 0.4, draw= black, double] {};
    }
\draw[dotted] (-1.5,0)  -- (0,0) node[pos=-0.5] {$t=0$};
\draw[dotted] (-1.5,5)  -- (0,5) node[pos=-0.5] {$t=L$};
\draw[dotted] (-1.5,10)  -- (0,10) node[pos=-0.5] {$t=2L$};
\draw[dotted] (0,-1)  -- (0,0) node[pos=-0.3] {$x=0$};
\draw[dotted] (5,-1)  -- (5,0) node[pos=-0.3] {$x=L$};

\end{tikzpicture}%
    \caption{Example mesh satisfying \eqref{QhVD to QhVOD}. On each time-slice, 
    the degrees of freedom (DoFs) of each boundary agree with the DoFs on the opposite boundary shifted in time by $L$.}
    \label{fig:Example_mesh}
\end{figure}

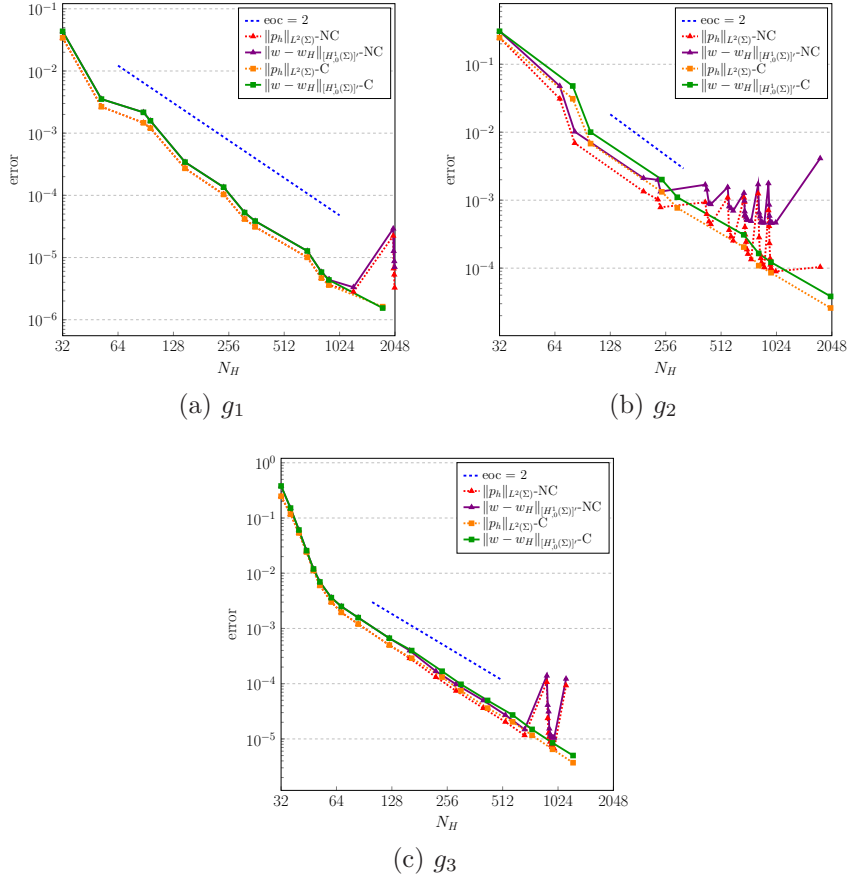
\begin{figure}[H]%
\vspace{-0.5cm}
 \centering
 \subfloat[$g_1$]{\pgfplotstabletranspose\noCon{Data/default_adapt_A_m3.txt}
\pgfplotstabletranspose\Con{Data/default_adapt_A_m3_Con.txt}


{\scalefont{1.5}
\begin{tikzpicture}
\begin{loglogaxis}[
    width=14cm,
    height=14cm,
    xlabel={$N_H$},
    ylabel={error},
    xmin=32, xmax=2048,
    ytick pos= left,
    xtick pos= bottom,
    xtick={32,64,128,256,512,1024,2048},
    xticklabels= {32,64,128,256,512,1024,2048},
    legend pos=north east,
    legend cell align={left},
    legend style ={font = \large, inner sep= 6pt, outer sep = -6pt},
    ymajorgrids=true,
    grid style=dashed,
]

\addplot [
    domain= 64:1024, 
    samples=8, 
    color=blue,
    dashed,
    line width=2pt
]
{50*x^-2};
\addlegendentry{eoc = 2}




\addplot[red,mark= triangle*,line width=2pt, mark options={solid}, dotted] table [x index=3, y index = 1] {\noCon};
\addlegendentry{$\|p_h\|_{L^2(\Sigma)}$-NC}

\addplot[violet, mark=triangle*,mark options=solid, line width=2pt] table [x index= 3, y index = 2] {\noCon};
\addlegendentry{$\|w-w_H\|_{[H_{,0}^1(\Sigma)]^\prime}$-NC}

\addplot[orange,mark= square*,line width=2pt, mark options={solid}, dotted] table [x index=3, y index = 1] {\Con};
\addlegendentry{$\|p_h\|_{L^2(\Sigma)}$-C}

\addplot[green!60!black, mark=square*,mark options=solid, line width=2pt] table [x index= 3, y index = 2] {\Con};
\addlegendentry{$\|w-w_H\|_{[H_{,0}^1(\Sigma)]^\prime}$-C}

\end{loglogaxis}
\end{tikzpicture}
} }%
 \subfloat[$g_2$]{\pgfplotstabletranspose\noCon{Data/default_adapt_B_m3.txt}
\pgfplotstabletranspose\Con{Data/default_adapt_B_m3_Con.txt}


{\scalefont{1.5}
\begin{tikzpicture}
\begin{loglogaxis}[
    width=14cm,
    height=14cm,
    xlabel={$N_H$},
    ylabel={error},
    xmin=32, xmax=2048,
    ytick pos= left,
    xtick pos= bottom,
    xtick={32,64,128,256,512,1024,2048},
    xticklabels= {32,64,128,256,512,1024,2048},
    legend pos=north east,
    legend cell align={left},
    legend style ={font = \large, inner sep= 6pt, outer sep = -6pt},
    ymajorgrids=true,
    grid style=dashed,
]

\addplot [
    domain= 128:320, 
    samples=8, 
    color=blue,
    dashed,
    line width=2pt
]
{300*x^-2};
\addlegendentry{eoc = 2}




\addplot[red,mark= triangle*,line width=2pt, mark options={solid}, dotted] table [x index=3, y index = 1] {\noCon};
\addlegendentry{$\|p_h\|_{L^2(\Sigma)}$-NC}

\addplot[violet, mark=triangle*,mark options=solid, line width=2pt] table [x index= 3, y index = 2] {\noCon};
\addlegendentry{$\|w-w_H\|_{[H_{,0}^1(\Sigma)]^\prime}$-NC}

\addplot[orange,mark= square*,line width=2pt, mark options={solid}, dotted] table [x index=3, y index = 1] {\Con};
\addlegendentry{$\|p_h\|_{L^2(\Sigma)}$-C}

\addplot[green!60!black, mark=square*,mark options=solid, line width=2pt] table [x index= 3, y index = 2] {\Con};
\addlegendentry{$\|w-w_H\|_{[H_{,0}^1(\Sigma)]^\prime}$-C}

\end{loglogaxis}
\end{tikzpicture}
}}\\
 \subfloat[$g_3$]{\pgfplotstabletranspose\noCon{Data/default_adapt_C_m3.txt}
\pgfplotstabletranspose\Con{Data/default_adapt_C_m3_Con.txt}


{\scalefont{1.5}
\begin{tikzpicture}
\begin{loglogaxis}[
    width=14cm,
    height=14cm,
    xlabel={$N_H$},
    ylabel={error},
    xmin=32, xmax=2048,
    ytick pos= left,
    xtick pos= bottom,
    xtick={32,64,128,256,512,1024,2048},
    xticklabels= {32,64,128,256,512,1024,2048},
    legend pos=north east,
    legend cell align={left},
    legend style ={font = \large, inner sep= 6pt, outer sep = -6pt},
    ymajorgrids=true,
    grid style=dashed,
]

\addplot [
    domain= 100:512, 
    samples=8, 
    color=blue,
    dashed,
    line width=2pt
]
{30*x^-2};
\addlegendentry{eoc = 2}




\addplot[red,mark= triangle*,line width=2pt, mark options={solid}, dotted] table [x index=3, y index = 1] {\noCon};
\addlegendentry{$\|p_h\|_{L^2(\Sigma)}$-NC}

\addplot[violet, mark=triangle*,mark options=solid, line width=2pt] table [x index= 3, y index = 2] {\noCon};
\addlegendentry{$\|w-w_H\|_{[H_{,0}^1(\Sigma)]^\prime}$-NC}

\addplot[orange,mark= square*,line width=2pt, mark options={solid}, dotted] table [x index=3, y index = 1] {\Con};
\addlegendentry{$\|p_h\|_{L^2(\Sigma)}$-C}

\addplot[green!60!black, mark=square*,mark options=solid, line width=2pt] table [x index= 3, y index = 2] {\Con};
\addlegendentry{$\|w-w_H\|_{[H_{,0}^1(\Sigma)]^\prime}$-C}

\end{loglogaxis}
\end{tikzpicture}
}}%
 \caption{Error convergence comparison between non-constrained (NC) and constrained(C) adaptive refinement algorithms for \eqref{vf:default}, given $m=3$ and different Dirichlet data.}%
 \vspace{-0.5cm}
 \label{fig:constraint_comparison}%
\end{figure}

The next point in our 'numerical agenda' is to study the need for
\emph{uniform meshes} in Theorem~\ref{thm:discrete-stability-n-time-slices}.
We remind the reader that discrete inf-sup stability of $\OV$ relies on
the assumption that \eqref{QhVD to QhVOD} is satisfied. As stated in
Remark \ref{rem:meshconstraints}, the adaptive procedure should uphold
the constraint on the mesh given by \eqref{QhVD to QhVOD} when the
initial mesh is uniform. In practice, however, it seems that at some
point during the refinement routine the constraint on the mesh is no
longer satisfied, resulting in a loss of discrete inf-sup stability.
This explains the inability of the adaptively refined formulation to
converge after a certain number of refinements, as visualised in
Fig.~\ref{Fig:natural_m3_AB}. To circumvent this issue, we consider a
\emph{constrained} adaptive algorithm with the mesh condition
\eqref{QhVD to QhVOD} hard-coded into the implementation. Ensuring that
the mesh on the boundary restricted to a time-slice corresponds to the mesh
on the opposite boundary of the subsequent time-slice, provides a sufficient
condition to satisfy \eqref{QhVD to QhVOD}. An example of a mesh satisfying
this condition is given in Fig.~\ref{fig:Example_mesh}. The constrained
adaptive refinement routine is realised by enforcing this condition at
each iteration. A comparison of the non-constrained and constrained adaptive
refinement routines is given in Fig.~\ref{fig:constraint_comparison}.
During the early stages of the refinement procedure, non-constrained
adaptive refinement may result in a higher convergence rate compared to
the constrained algorithm. This can be explained by the fact the constrained
refinement scheme may unnecessarily refine parts where the Galerkin solution
is zero. After several refinements, the constrained algorithm overcomes
the issue encountered by its unconstrained counterpart.

Finally, we compare the performance of the proposed adaptive
algorithm from formulation \eqref{vf:EBEM}, which we will dub \emph{LSBEM}, with an adaptive BEM routine
introduced in \cite{SchulzSteinbach:2000} and applied to the wave equation
in \cite{SteinbachZank:2016}, here referred to as \emph{SteZan}. Performance
is measured by considering the error with respect to the amount of degrees
of freedom. For the numerical experiments a Galerkin approximation of the
\emph{direct} energetic BIE is considered:
For $i\in\{1,2\}$, find $z_H \in \pwc_H^0(\Sigma) \subset L^2(\Sigma)$ such that
\begin{equation}\label{eq:directEBEM}
  \langle\partial_t \OV w_H,q_H\rangle_\Sigma =
  \frac{1}{2}\langle\partial_t g_i,q_H\rangle_\Sigma +
  \langle\partial_t \OK g_i,q_H\rangle_\Sigma, \quad
  \forall q_H \in \pwc_H^0(\Sigma)\subset L^2(\Sigma),
\end{equation}
where $\OK$ denotes the \emph{double layer operator}, which is given for
$g=0$ outside of $\Sigma$ by \cite{SteinbachZank:2016}:
\begin{equation*}\label{def:K}
\OK g(x,t) = \begin{cases}
              -\frac{1}{2}g(L,t-L) & x=0,\\[1mm]
              -\frac{1}{2}g(0,t-L) & x=L.
             \end{cases}
\end{equation*}

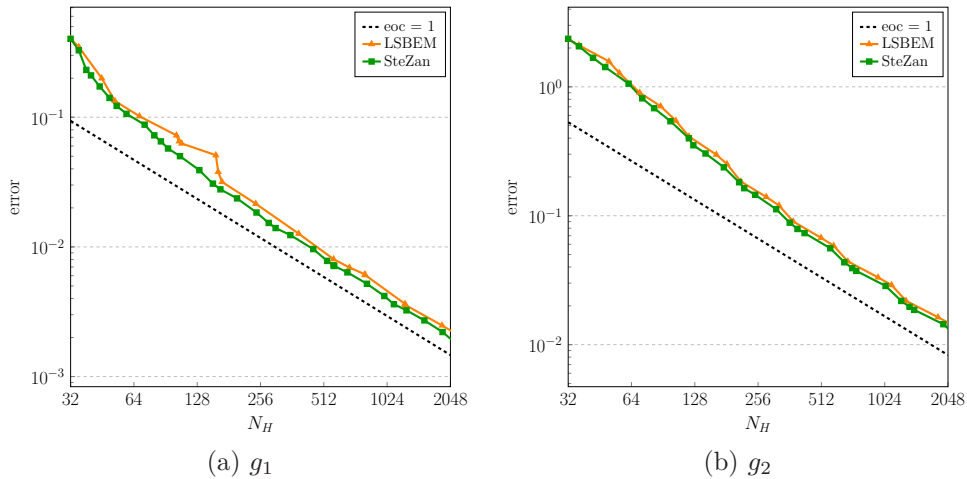
\begin{figure}[H]
    \centering
    \subfloat[][\centering $g_1$]{{\pgfplotstabletranspose\LSBEM{Data/EBEM_adapt_A_Direct_m2.txt}
\pgfplotstabletranspose\SteZan{Data/EBEM_adaptv2_A_Direct_m1.txt}


{\scalefont{1.5}
\begin{tikzpicture}
\begin{loglogaxis}[
    width=14cm,
    height=14cm,
    xlabel={$N_H$},
    ylabel={error},
    xmin=32, xmax=2048,
    ytick pos= left,
    xtick pos= bottom,
    xtick={32,64,128,256,512,1024,2048},
    xticklabels= {32,64,128,256,512,1024,2048},
    legend pos=north east,
    legend cell align={left},
    legend style ={font = \large, inner sep= 6pt, outer sep = -6pt},
    ymajorgrids=true,
    grid style=dashed,
]


\addplot [
    domain= 32:2048, 
    samples=8, 
    color=black,
    dashed,
    line width=2pt
]
{3*x^-1};
\addlegendentry{eoc = 1}





\addplot[orange, mark=triangle*,mark options=solid, line width=2pt] table [x index= 3, y index = 2] {\LSBEM};
\addlegendentry{LSBEM}


\addplot[green!60!black, mark=square*,mark options=solid, line width=2pt] table [x index= 3, y index = 2] {\SteZan};
\addlegendentry{SteZan}

\end{loglogaxis}
\end{tikzpicture}
} }}%
    \subfloat[][\centering $g_2$]{{\pgfplotstabletranspose\LSBEM{Data/EBEM_adapt_B_Direct_m2.txt}
\pgfplotstabletranspose\SteZan{Data/EBEM_adaptv2_B_Direct_m1.txt}


{\scalefont{1.5}
\begin{tikzpicture}
\begin{loglogaxis}[
    width=14cm,
    height=14cm,
    xlabel={$N_H$},
    ylabel={error},
    xmin=32, xmax=2048,
    ytick pos= left,
    xtick pos= bottom,
    xtick={32,64,128,256,512,1024,2048},
    xticklabels= {32,64,128,256,512,1024,2048},
    legend pos=north east,
    legend cell align={left},
    legend style ={font = \large, inner sep= 6pt, outer sep = -6pt},
    ymajorgrids=true,
    grid style=dashed,
]


\addplot [
    domain= 32:2048, 
    samples=8, 
    color=black,
    dashed,
    line width=2pt
]
{17*x^-1};
\addlegendentry{eoc = 1}





\addplot[orange, mark=triangle*,mark options=solid, line width=2pt] table [x index= 3, y index = 2] {\LSBEM};
\addlegendentry{LSBEM}


\addplot[green!60!black, mark=square*,mark options=solid, line width=2pt] table [x index= 3, y index = 2] {\SteZan};
\addlegendentry{SteZan}

\end{loglogaxis}
\end{tikzpicture}
} }}%
    \caption{ $L^2(\Sigma)$-error convergence comparison for proposed adaptive algorithm (LSBEM) and method from the literature (SteZan). }
    \label{fig:EBEM_adaptive_comparison}
\end{figure}

A comparison of both methods is presented in
Fig.~\ref{fig:EBEM_adaptive_comparison} using Dirichlet data $g_1$ and
$g_2$. Both methods perform similarly when it comes to error convergence.
However, in order to obtain a valid error estimator, the LSBEM approach
requires solving for a mixed boundary element method with $m=2$, increasing
the computational complexity at each refinement. On the other hand, the
SteZan method has limited applicability: it is restricted to direct
formulations and requires an implementation or approximation of the
adjoint double layer and hypersingular operators.

\subsection{Stability Constant}\label{ssec: stability constant}
\begin{figure}[H]
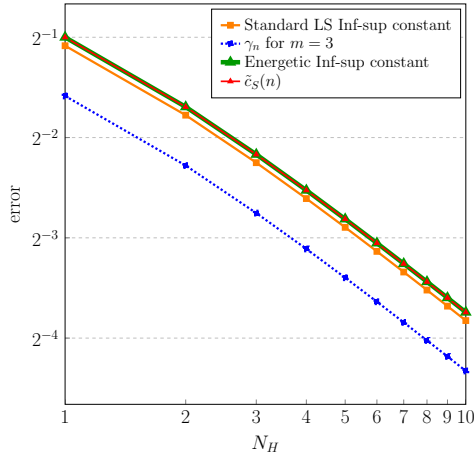

    \centering
    \includestandalone[width=0.4\textwidth]{Figures/inf_sup_constant}%
    \caption{Comparison between $\gamma_n$ and $\tilde c_S(n)$, and the discrete inf-sup constant for formulations \eqref{vf:default} with $m=3$ and \eqref{vf:EBEM} with $m=1$, on different amount of time-slices.}
    \label{fig:inf_sup_default}
\end{figure}

Finally, we compare the theoretical stability constant, as proposed in
Theorem \ref{thm:discrete-stability-n-time-slices}, with the actual discrete
inf-sup constant, computed using the method introduced in, e.g.,
\cite[Rem.~3.159]{John:2016}. For notational convenience let us denote the
theoretical stability constant by
\begin{equation*}\label{defn: notation stability constant}
  \gamma_{n}:= 2\sin\left(\frac{\pi}{2(2n+1)}\right)
  \left(\frac{1}{2}-\frac{1}{m}\right).
\end{equation*}
The computation of the discrete inf-sup constant requires the usage of a
mass matrix with respect to the $[H^1_{,0}(\Sigma)]^\prime$-inner product,
the implementation of this matrix is based on the theory presented in
Section \ref{ssec: dual norm}. The results for the stable ($m=3$) standard
formulation \eqref{vf:default} and the energetic BEM formulation
\eqref{vf:EBEM} without nesting $(m=1)$ are given in
Fig.~\ref{fig:inf_sup_default}. There we see that the proposed stability
constant has the same asymptotic behaviour as the actual discrete inf-sup
constant. In the case of energetic BEM, which is stable for $m=1$, we
observe that the stability constant coincides with $\tilde{c}_S(n)$, as
defined in Theorem \ref{Thm ellipticity V*V}. On each time-slice the
coarse mesh consists of $32$ uniform elements.

\section{Conclusions}
In this paper we have formulated and analyzed a least squares approach
for first kind boundary integral equations for the Dirichlet problem
for the wave equation. We have established stability of a related
boundary element method, from which we can derive a priori error estimates.
Moreover, the approximation of the adjoint variable can be used as
an error indicator to drive an adaptive algorithm. Numerical results,
also for less regular Dirichlet data, confirm the theoretical findings.

It is more or less obvious that this approach can be applied as well to
problems with different boundary conditions, and to other boundary
integral equations also including the double layer operator and its
adjoint, and the hypersingular boundary integral operator for the
wave equation. A possible extension to systems such as in elastodynamics
will also follow the lines as given for the scalar wave equation.
More challenging is the construction of efficient solution methods for the
resulting linear systems of algebraic equations, and the construction
of appropriate preconditioners. The implementation of the proposed
approach to solve problems in higher space dimensions is ongoing work,
but the numerical analysis can not be done in such an explicit way as
it is possible in one dimension.


\end{document}


\pgfplotstableread{Data/inf_sup_default.txt}\infsupDefault
\pgfplotstableread{Data/inf_sup_EBEM.txt}\infsupEBEM


{\scalefont{1.5}
\begin{tikzpicture}
\begin{loglogaxis}[
    log basis y ={2},
    log basis x = {2},
    width=14cm,
    height=14cm,
    xlabel={$N_H$},
    ylabel={error},
    xmin=1, xmax=10,
    ytick pos= left,
    xtick pos= bottom,
    xtick={1,2,3,4,5,6,7,8,9,10},
    xticklabels= {1,2,3,4,5,6,7,8,9,10},
    ytick = {0.5,0.25,0.125, 0.0625},
    legend pos=north east,
    legend cell align={left},
    legend style ={font = \large, inner sep= 6pt, outer sep = -6pt},
    ymajorgrids=true,
    grid style=dashed,
]

\addplot[orange,mark= square*,line width=2pt, mark options={solid}] table [x expr=\coordindex+1, y index = 0] {\infsupDefault};
\addlegendentry{Standard LS Inf-sup constant}

\addplot [
domain= 1:10, 
samples at={0,1,2,3,4,5,6,7,8,9,10}, 
color=blue,
dotted,
mark=square*,
line width=2pt
]
{(1/2-1/6)*2*sin(deg(pi/(4*x+2)))};
\addlegendentry{$\gamma_n$ for $m=3$}

\addplot[green!60!black, mark=triangle*,mark options=solid, line width=4pt] table [x expr = \coordindex+1, y index = 0] {\infsupEBEM};
\addlegendentry{Energetic Inf-sup constant}

\addplot [
domain= 1:10, 
samples at={0,1,2,3,4,5,6,7,8,9,10}, 
color=red,
mark=triangle*,
line width=1.5pt
]
{sin(deg(pi/(4*x+2)))};
\addlegendentry{$\tilde c_S(n)$}




\end{loglogaxis}
\end{tikzpicture}
}